\documentclass[12pt]{amsart}
\usepackage{amsmath}
\usepackage{enumerate,amssymb,amsfonts,latexsym}
\usepackage{amscd, mathtools, tabu}
\usepackage{paralist}
\usepackage[misc]{ifsym}
\usepackage{epsfig}
\usepackage{epstopdf}
\usepackage[unicode, pdfproducer={Latex with hyperref}, pdfcreator={pdflatex}, colorlinks=true]{hyperref}
\hypersetup{urlcolor=blue, citecolor=red}
\usepackage{graphicx}
\usepackage{xparse}
\usepackage[utf8]{inputenc}
\usepackage[T1]{fontenc}
\usepackage{lmodern}
\usepackage{float}
\usepackage{tikz}
\usepackage{mathtools}
\usepackage{csquotes}
\usepackage{scalerel,stackengine}
\allowdisplaybreaks

\newtheorem{theorem}{Theorem}[section]
\newtheorem{corollary}[theorem]{Corollary}

\newtheorem{lemma}[theorem]{Lemma}
\newtheorem{proposition}[theorem]{Proposition}

\theoremstyle{definition}
\newtheorem{definition}[theorem]{Definition}
\newtheorem{remark}[theorem]{Remark}

\title[Geometry of area-preserving maps]
{On the geometry of period doubling invariant sets for area-preserving maps}

\author[Denis Gaidashev and Dan Lilja]{}

\thanks{The first author is supported by Vetenskapsr\r{a}det project grant 2022-05046}
\thanks{$^*$Corresponding author: Denis Gaidashev}

\begin{document}
\maketitle

\centerline{\scshape Denis Gaidashev$^{{\href{mailto:gaidash@math.uu.se}{\textrm{\Letter}}}*1}$ and Dan Lilja$^{{\href{mailto:dan@danlilja.se}{\textrm{\Letter}}}2}$}

\medskip

{\footnotesize
	\centerline{$^1$Department of Mathematics, Uppsala University, Uppsala, Sweden}
}

\medskip

{\footnotesize
	\centerline{$^2$Consid, Sweden}
}

\bigskip

\begin{abstract} 
	The geometry of the period doubling Cantor sets of strongly dissipative
	infinitely renormalizable Hénon-like maps has been shown to be unbounded by M.
	Lyubich, M. Martens and A. de Carvalho, although the measure of unbounded
	``spots'' in the Cantor set has been demonstrated to be zero.

	We show that an even more extreme situation takes places for infinitely
	renormalizable area-preserving Hénon-like maps: both bounded and unbounded
	geometries exist on subsets of positive measure.
\end{abstract}

\newtheorem{thm}{Theorem}[section]
\newtheorem*{thmA}{Theorem A}
\newtheorem*{thmB}{Theorem B}
\newtheorem*{thmC}{Theorem C}
\newtheorem{propZA}{Proposition: Zero Cone Expansion on the Average}
\newtheorem{clm}[thm]{Claim}

\newtheorem*{zeroconeexpansion}{Proposition: Zero Cone Expansion on the Average}

\numberwithin{equation}{section}
\numberwithin{figure}{section}

\newcommand{\lab}[3]{\psfrag{#1}[#3]{$\scriptstyle{#2}$}}
\newcommand{\lowlimit}[1]{\substack{\phantom{-} \\ #1}}

\font\nt=cmr7

\def\note#1
{\marginpar
{\nt $\leftarrow$
\par
\hfuzz=20pt \hbadness=9000 \hyphenpenalty=-100 \exhyphenpenalty=-100
\pretolerance=-1 \tolerance=9999 \doublehyphendemerits=-100000
\finalhyphendemerits=-100000 \baselineskip=6pt
#1}\hfuzz=1pt}

\newcommand{\bignote}[1]{\begin{quote} \sf #1 \end{quote}}

\newcommand{\vertiii}[1]{{\left\vert\kern-0.25ex\left\vert\kern-0.25ex\left\vert #1 
    \right\vert\kern-0.25ex\right\vert\kern-0.25ex\right\vert}}

\makeatletter
\newcommand{\raisemath}[1]{\mathpalette{\raisem@th{#1}}}
\newcommand{\raisem@th}[3]{\raisebox{#1}{$#2#3$}}
\makeatother
\newcommand{\move}[1]{\hspace{-0.4mm}\raisemath{2pt}{#1} }
\newcommand{\touch}[1]{\hspace{-0.2mm}{\raisebox{0.4pt}{$\scriptstyle{#1}$}} }

\stackMath
\newcommand\hhat[1]{%
\savestack{\tmpbox}{\stretchto{%
  \scaleto{%
    \scalerel*[\widthof{\ensuremath{#1}}]{\kern-.6pt\bigwedge\kern-.6pt}%
    {\rule[-\textheight/2]{1ex}{\textheight}}
  }{\textheight}%
}{0.5ex}}%
\stackon[1pt]{#1}{\tmpbox}%
}
\parskip 1ex

\newcommand{\QED}{\rlap{$\sqcup$}$\sqcap$\smallskip}

\renewcommand{\Im}{\operatorname{Im}}

\newcommand{\shortplus}{\scalebox{0.5}[0.5]{\( + \)}}
\newcommand{\shortminus}{\scalebox{0.75}[1.0]{\( - \)}}

\def\sss{\subsubsection}

\newcommand{\di}{\partial}
\newcommand{\dibar}{\bar\partial}
\newcommand{\hookra}{\hookrightarrow}
\newcommand{\ra}{\rightarrow}
\newcommand{\hra}{\hookrightarrow}
\newcommand{\imply}{\Rightarrow}
\def\lra{\longrightarrow}
\newcommand{\wc}{\underset{w}{\to}}
\newcommand{\tu}{\textup}

\def\ssk{\smallskip}
\def\msk{\medskip}
\def\bsk{\bigskip}
\def\noi{\noindent}
\def\nin{\noindent}
\def\lqq{\lq\lq}
\def\sm{\smallsetminus}
\def\bolshe{\succ}
\def\ssm{\smallsetminus}
\def\tr{{\text{tr}}}

\newcommand{\diam}{\operatorname{diam}}
\newcommand{\dist}{\operatorname{dist}}
\newcommand{\Hdist}{\operatorname{H-dist}}
\newcommand{\cl}{\operatorname{cl}}
\newcommand{\inter}{\operatorname{int}}
\renewcommand{\mod}{\operatorname{mod}}
\newcommand{\card}{\operatorname{card}}
\newcommand{\tl}{\tilde}
\newcommand{\ind}{ \operatorname{ind} }
\newcommand{\Dist}{\operatorname{Dist}}
\newcommand{\len}{\operatorname{\l}}

\newcommand{\ctg}{\operatorname{ctg}}
\newcommand{\arcctg}{\operatorname{arcctg}}
\newcommand{\orb}{\operatorname{orb}}
\newcommand{\HD}{\operatorname{HD}}
\newcommand{\supp}{\operatorname{supp}}
\newcommand{\id}{\operatorname{id}}
\newcommand{\length}{\operatorname{length}}
\newcommand{\dens}{\operatorname{dens}}
\newcommand{\meas}{\operatorname{meas}}
\newcommand{\distM}{\operatorname{dist}_{Mon}} 
\newcommand{\per}{\operatorname{per}}

\renewcommand{\d}{{\diamond}}

\newcommand{\Dil}{\operatorname{Dil}}
\newcommand{\Ker}{\operatorname{Ker}}
\newcommand{\tg}{\operatorname{tg}}
\newcommand{\codim}{\operatorname{codim}}
\newcommand{\isom}{\approx}
\newcommand{\comp}{\circ}
\newcommand{\esssup}{\operatorname{ess-sup}}

\newcommand{\transverse}
 {\kern .7em\makebox[0pt][c]{\raisebox{.2ex}{$|$}}\kern -.6em\cap}

\newcommand{\tangent}
 {\kern .7em\makebox[0pt][c]{\raisebox{.77ex}{$--$}}\kern -.6em\cap}

\newcommand{\SLa}{\underset{\La}{\Subset}}

\newcommand{\const}{\mathrm{const}}
\def\loc{{\mathrm{loc}}}
\def\fib{{\mathrm{fib}}}

\newcommand{\bm}{\bar{m}}
\newcommand{\eps}{{\varepsilon}}
\newcommand{\epsi}{{\epsilon}}
\newcommand{\veps}{{\varepsilon}}
\newcommand{\De}{{\Delta}}
\newcommand{\de}{{\delta}}
\newcommand{\la}{{\lambda}}
\newcommand{\La}{{\Lambda}}
\newcommand{\si}{{\sigma}}
\newcommand{\Si}{{\Sigma}}
\newcommand{\Om}{{\Omega}}
\newcommand{\om}{{\omega}}
\newcommand{\Wl}{W^s_{\mathrm{loc}}(F_*)}
\newcommand{\Wll}{W^{ss}_{\matrm{loc}}(F_*)}

\newcommand{\al}{{\alpha}}
\newcommand{\ba}{{\mbox{\boldmath$\alpha$} }}
\newcommand{\bbe}{{\mbox{\boldmath$\beta$} }}
\newcommand{\bk}{{\boldsymbol{\kappa}}}
\newcommand{\bg}{{\boldsymbol{\gamma}}}

\newcommand{\bare}{{\bar\eps}}

\newcommand{\Ray}{{\mathcal R}}
\newcommand{\Eq}{{\mathcal E}}
\newcommand{\PR}{PR}

\newcommand{\AAA}{{\mathcal A}}
\newcommand{\BB}{{\mathcal B}}
\newcommand{\CC}{{\mathcal C}}
\newcommand{\DD}{{\mathcal D}}
\newcommand{\EE}{{\mathcal E}}
\newcommand{\EEE}{{\mathcal O}}
\newcommand{\II}{{\mathcal I}}
\newcommand{\FF}{{\mathcal F}}
\newcommand{\GG}{{\mathcal G}}
\newcommand{\JJ}{{\mathcal J}}
\newcommand{\HH}{{\mathcal H}}
\newcommand{\KK}{{\mathcal K}}
\newcommand{\LL}{{\mathcal L}}
\newcommand{\MM}{{\mathcal M}}
\newcommand{\NN}{{\mathcal N}}
\newcommand{\OO}{{\mathcal O}}
\newcommand{\PP}{{\mathcal P}}
\newcommand{\QQ}{{\mathcal Q}}
\newcommand{\QM}{{\mathcal QM}}
\newcommand{\QP}{{\mathcal QP}}
\newcommand{\QL}{{\mathcal QL}}

\newcommand{\RR}{{\mathcal R}}
\renewcommand{\SS}{{\mathcal S}}
\newcommand{\TT}{{\mathcal T}}
\newcommand{\TTT}{{\mathcal P}}
\newcommand{\UU}{{\mathcal U}}
\newcommand{\VV}{{\mathcal V}}
\newcommand{\WW}{{\mathcal W}}
\newcommand{\XX}{{\mathcal X}}
\newcommand{\YY}{{\mathcal Y}}
\newcommand{\ZZ}{{\mathcal Z}}

\newcommand{\A}{{\mathbb A}}
\newcommand{\E}{{\mathbb E}}    
\newcommand{\C}{{\mathbb C}}
\newcommand{\bC}{{\bar{\mathbb C}}}
\newcommand{\D}{{\mathbb D}}
\newcommand{\Hyp}{{\mathbb H}}
\newcommand{\J}{{\mathbb J}}
\newcommand{\I}{{\mathbb I}}    
\newcommand{\Ll}{{\mathbb L}}
\renewcommand{\L}{{\mathbb L}}
\newcommand{\M}{{\mathbb M}}
\newcommand{\N}{{\mathbb N}}
\newcommand{\Q}{{\mathbb Q}}
\newcommand{\R}{{\mathbb R}}
\newcommand{\T}{{\mathbb T}}
\newcommand{\V}{{\mathbb V}}
\newcommand{\U}{{\mathbb U}}
\newcommand{\W}{{\mathbb W}}
\newcommand{\X}{{\mathbb X}}
\newcommand{\Z}{{\mathbb Z}}

\newcommand{\VVV}{{\mathbf U}}
\newcommand{\UUU}{{\mathbf U}}

\newcommand{\tT}{{\mathrm{T}}}
\newcommand{\tD}{{D}}
\newcommand{\hyp}{{\mathrm{hyp}}}

\newcommand{\f}{{\bf f}}
\newcommand{\g}{{\bf g}}
\newcommand{\bh}{{\bf h}}
\renewcommand{\i}{{\bar i}}
\renewcommand{\j}{{\bar j}}
\newcommand{\bv} {{\bf v  }}
\newcommand{\bw} {{\bf w  }}
\renewcommand{\k}{\kappa}

\def\Bf{{\mathbf{f}}}
\def\Bg{{\mathbf{g}}}
\def\BG{{\mathbf{G}}}
\def\Bh{{\mathbf{h}}}
\def\Bv{{\mathbf{v}}}
\def\Bw{{\mathbf{w}}}
\def\Bz{{\mathbf{z}}}
\def\Bx{{\mathbf{x}}}
\def\By{{\mathbf{y}}}

\def\Bde{{\boldsymbol{\de}}}

\def\bB{{\mathbf{B}}}
\def\bQ{{\mathbf{Q}}}
\def\BH{{\mathbf{H}}}
\def\BF{{\mathbf{F}}}
\def\BS{{\mathbf{S}}}
\def\BDe{{\boldsymbol{\Delta}}}
\def\BGa{{\boldsymbol{\Gamma}}}

\def\BT{{\mathbf{T}}}
\def\Bj{{\mathbf{j}}}
\def\Bphi{{\mathbf{\Phi}}}
\def\BPsi{{\boldsymbol{\Psi}}}
\def\BPhi{{\boldsymbol{\BPhi}}}
\def\B0{{\mathbf{0}}}
\def\BU{{\mathbf{U}}}
\def\BV{{\mathbf{V}}}
\def\BR{{\mathbf{R}}}
\def\BG{{\mathbf{G}}}
\newcommand{\Comb}{{\it Comb}}
\newcommand{\Top}{{\it Top}}
\newcommand{\QC}{\mathcal QC}
\newcommand{\Def}{\mathcal Def}
\newcommand{\Teich}{\mathcal Teich}
\newcommand{\PPL}{{\mathcal P}{\mathcal L}}
\newcommand{\Jac}{\operatorname{Jac}}
\renewcommand{\DH}{\operatorname{DH}}
\newcommand{\Homeo}{\operatorname{Homeo}}
\newcommand{\AC}{\operatorname{AC}}
\newcommand{\Dom}{\operatorname{Dom}}

\newcommand{\Aff}{\operatorname{Aff}}
\newcommand{\Euc}{\operatorname{Euc}}
\newcommand{\MobC}{\operatorname{M\ddot{o}b}({\mathbb C}) }
\newcommand{\PSL}{ {\mathcal{PSL}} }
\newcommand{\SL}{ {\mathcal{SL}} }
\newcommand{\CP}{ {\mathbb{CP}}   }

\newcommand{\hf}{{\hat f}}
\newcommand{\hz}{{\hat z}}
\newcommand{\hM}{{\hat M}}

\renewcommand{\lq}{``}
\renewcommand{\rq}{''}

\newcommand{\Domain}{\operatorname{Dom}}
\newcommand{\Image}{\operatorname{Im}}
\newcommand{\graph}{\operatorname{graph}}
\newcommand{\Trap}{\operatorname{Trap}}
\newcommand{\Orb}{\operatorname{Orb}}

\def\diff{{\text{diff}}}
\def\aff{{\text{aff}}}


\catcode`\@=12

\def\Empty{}
\newcommand\oplabel[1]{
  \def\OpArg{#1} \ifx \OpArg\Empty {} \else
  	\label{#1}
  \fi}
		
%


%

\newcommand{\comm}[1]{}
\newcommand{\comment}[1]{}

\newcommand{\mm}{\hspace{-1.1mm}-\hspace{-1.1mm}}
\newcommand{\pp}{\hspace{-1.1mm}+\hspace{-1.1mm}}       

\NewDocumentCommand{\ceil}{s O{} m}{%
  \IfBooleanTF{#1} 
    {\left\lceil#3\right\rceil} 
    {#2\lceil#3#2\rceil} 
}


\setcounter{tocdepth}{1}


\section{Introduction}\label{intro}

The period doubling cascade is one of the fundamental scenarios of transition
from periodic to chaotic dynamics in one-dimensional systems. This cascade
accumulates on a dynamical system that admits an attracting invariant Cantor
set. The properties of these period doubling Cantor sets are very well
understood. In the late 1970's Feigenbaum~\cite{F1,F2} and, independently,
Coullet and Tresser~\cite{CT,TC}, discovered numerically universal geometric
properties of these Cantor sets.
 
M. Feigenbaum, P. Coullet and C. Tresser introduced renormalization in dynamics
to explain the observed geometrical universality. Renormalization, viewed as an
operator on a class of dynamical systems, maps one such system into another one
which corresponds to a rescaled version of a higher iteration of the original
system acting on a subset of its phase space. This renormalization operator
typically has a hyperbolic horseshoe; the dynamics of the systems on the stable
leaves converges to an orbit (which can be periodic), and the behavior of the
renormalization operator around this orbit determines the asymptotic small scale
properties. This explains the observed universality.

The renormalization technique has been generalized to many types of dynamics.
However, a rigorous study of universality has been surprisingly difficult and
technically sophisticated. It has been thoroughly carried out in the case
of one-dimensional maps, on the interval or the circle, see~\cite{AL,  FMP,He,L,Ma,McM,MS, S,VSK,Y}.

One of the strong properties of infinitely renormalizable maps is
\textit{rigidity}: asymptotically, on smaller and smaller scales, there is a
universal geometry around the invariant Cantor set. All period doubling Cantor
sets are topologically equivalent, but a priori there is no reason to believe
that these conjugating homeomorphisms can have some smoothness. However, it is
well-known that the period doubling Cantor sets in one-dimensional dynamics are
rigid: there are smooth conjugations.

Many numerical and physical experiments show that exactly the same universal
geometry from one-dimensional dynamics occurs also in some dissipative higher
dimensional systems. As it turns out, the rigidity phenomenon is more delicate
in higher dimensions.

Recently, one-dimensional techniques have been extended to strongly dissipative
perturbations of one-dimensional systems, such as Hénon maps, see~\cite{CEK1,
  CLM, LM1, LM2,Haz,HLM}. Strongly dissipative two-dimensional Hénon-like maps
can be thought of as two-dimensional perturbations of one-dimensional systems.
In~\cite{CEK1} and~\cite{CLM} two renormalization schemes were developed for
strongly dissipative Hénon-like maps at the accumulation of period doubling
which explain the universal geometry present in these Hénon-like maps.
Surprisingly, the period doubling Hénon-like Cantor sets are not smoothly
conjugated to their one-dimensional counterpart, see~\cite{CLM, LM1}.

Nevertheless, the geometry of the one-dimensional period doubling Cantor set is
still present. The conjugations between the Cantor sets of strongly dissipative
Hénon-like maps is almost everywhere, with respect to the natural measure on the
Cantor set, smooth. This phenomenon is called \textit{Probabilistic Rigidity},
see~\cite{LM1, LM3}. Small scale geometry has a probabilistic nature in higher
dimensions.

We will say that a Cantor set has bounded geometry, if the distance between
neighboring sets in the $n$-th generation cover is commensurable to their
diameters. These covers are generated dynamically by an iterated function system of two contractions, on of them being a constant linear map $\Lambda$, the other - a composition of this linear map with dynamics. Such iterated function system is commonly referred to as a {\it renormalization microscope}. Every  $n$-th generation cover contains a distinguished ``piece'', generated by a finite composition of $\Lambda$'s only. The $n$-th generation cover itself can be also viewed as the orbit of this central piece of length $2^n$ by dynamics.  

The papers~\cite{LM1} and~\cite{LM3} demonstrate that the unbounded
geometry is present in the Cantor sets of strongly dissipative Hénon-like maps,
but the measure of pieces in the $n$-th cover with unbounded geometry disappears
as $n$ increases.

The other extreme case is that of the area-preserving maps. Area-preserving maps
at accumulation of period doubling are observed by several authors in the early
80's, see~\cite{DP,Hl,BCGG,Bo,CEK2,EKW1}. In~\cite{EKW2} Eckmann, Koch and
Wittwer introduced a period doubling renormalization scheme for area preserving
maps and described the hyperbolic behavior of the renormalization operator in a
neighborhood of a renormalization fixed point. In particular, they observed
universality for maps at the accumulation of period doubling.

It was shown in~\cite{GJ1} that the maps in this Eckmann-Koch-Wittwer
universality class do have a period doubling Cantor set and the Lyapunov
exponents of dynamics restricted to this Cantor set are zero. It was later
demonstrated in~\cite{GJM} that the period doubling Cantor sets of the maps in
the stable manifold of the Eckmann-Koch-Wittwer renormalization fixed point are
{\it rigidly conjugate} by a $C^{1+\alpha}$ coordinate change.

We would like to emphasize that the proof of the rigidity of the period doubling Cantor set in ~\cite{GJM} uses in an essential way existence and hyperbolicity of the renormalization fixed point, demonstrated in ~\cite{EKW2} and ~\cite{GJ2}. Both of these works are computer-assisted. To date, there is no analytic proof of renormalization hyperbolicity in the area-preserving setting.

In our work we follow the approach of ~\cite{GJM} of providing analytic proofs which rely on previously obtained computer-assisted bounds.

In this paper we consider the geometry of the period doubling Cantor sets for
area-preserving infinitely renormalizable maps and demonstrate that the
situation is more complicated than for the dissipative ones. All $n$-th level
covers contain subsets of pieces with both bounded and unbounded
geometry, 
{ and the measure of each remains  positive as $n$ increases.}

\bigskip

\noindent {\bf Coexistence of Bounded and Unbounded Geometry.} {\it The period
  doubling Cantor sets of area-preserving maps in the Eckmann-Koch-Wittwer
  universality class have bounded and unbounded geometry on subsets of positive measure. }

\bigskip

The exact statements about the measure of bounded and unbounded pieces are
contained in Theorems A and B respectively.

We conclude this paper with a study of some ergodic properties of the derivative
cocycle of the fixed point map.

It turns out that the invariance of the fixed point map under the
renormalization operator, or, more generally, being in the renormalization
stable manifold, allows for a much better control of
convergence of ergodic means, or expansion/contraction within the pieces. We
will give informal statements of some of these consequences here, the exact
formulations can be found in the respective propositions.

\bigskip

\noindent {\bf Ergodic properties of the derivative cocycle.} Let $F_*$ denote the Eckmann-Koch-Wittwer fixed point.

\noindent {\it $(1)$ (\underline{Proposition $\ref{adequate_conv}$}.) For any
  function $f \in C \left(\R P^1 \times \{0,1\}^\N \right)$ the Birkhoff
  averages under the derivative cocycle over the odometer, $\Psi$, converge
  ``adequately'' to the space average $\bar f$, i.e. for a.a. $([u],w) \in \R
  P^1 \times \{0,1\}^\N)$,}
\begin{equation*}
  {1 \over 2^{n}}  \sum_{k=0}^{2^n-1}  f \circ \Psi^k([u],w) - \bar f=O \left(\beta^n n^{{3 \over 2}+\delta} \right)
\end{equation*}
{\it for some $0<\beta <1$ and any $\delta>0$.}

\bigskip

\noindent {\it $(2)$ (\underline{Proposition
    $\ref{zero_average_cone_expansion}$}.) The derivative cocycle preserves
  angles on the average: For a.a. $([u], w) \in \R P^1 \times \{0,1\}^\N$ and
  a.a. $([v], w) \in \R P^1 \times \{0,1\}^\N$,
  \[
    \lim_{k \rightarrow \infty} {1 \over k} \ln { \sin \left(D F^k_*(\tau_\rho) u, D F^k_*(\tau_\rho)  v \right) \over \sin(u, v)}=0,
  \]
  where $\tau_\rho$ is a point in the invariant Cantor set.}
\bigskip

\noindent {\it $(3)$ (\underline{Proposition $\ref{angles3}$}.) The distortion
  of the derivative $D F_*^k$ of long iterates evaluated at any two points $x$
  and $y$ of the ``central'' piece of the
  $n$-th cover of the Cantor set is bounded ``on the average'':
  \begin{equation} \label{avlogdist}
    {1 \over C}  < {1 \over 2^{n}} \sum_{k=0}^{2^n-1} \ln { \|  D F^k_*(x) u\| \over \|  D F^k_*(y) u\| } < C
  \end{equation}
  for some $C>0$ and a.a. non-zero vectors $u$.}

\bigskip

Since every piece in the $n$-th generation cover is an image of the central
piece by dynamics, result $(3)$, in particular, shows that ``on the average''
all pieces are bounded distortion diffeomorphic copies of the central piece.

\section{Preliminaries}\label{prelim}
Given a domain $\mathcal{D}\subset \mathbb{C}^2$, let $D\subset
\inter(\mathcal{D}\cap \mathbb{R}^2)$ be compactly contained in the real slice.
Assume $(0,0)\in D$. An \textit{area-preserving map} $F:D\to F(D)\subset
\mathbb{R}^2$ will mean a real analytic map which has a holomorphic extension to
$\mathcal{D}$, continuous on the boundary of $\DD$, and is an exact symplectic
diffeomorphism onto its image with the following properties

\medskip

\begin{itemize}
\item[1)] $T \circ F \circ T=F^{-1}$,  where $T(x,u)=(x,-u)$ (\textit{reversibility}),
\item[2)] $\partial_u y \ne 0$ with $(y,v)=F(x,u)$ (\textit{twist condition}).
\end{itemize}

\medskip

The collection of such maps is denoted by ${\rm Cons}(\DD)$. It has been shown
in~\cite{EKW2} that the set ${\rm Cons}(\DD)$ can be identified with a Banach
space $\AAA(\DD_s)$ of real symmetric functions $s: \DD_s \subset \C^2 \mapsto
\C$ holomorphic on some domain $\DD_s$, continuous on the boundary of $\DD_s$.
Specifically, $F \in {\rm Cons}(\DD)$ is generated by $s$:
\begin{equation}\label{sdef}
  \left({x  \atop  -s(y,x)} \right)  {{ \mbox{{\small{\textit{F}}}} \atop \mapsto} \atop \phantom{\mbox{\tiny{.}}}} \left({y \atop s(x,y) }\right).
\end{equation}

\subsection{Renormalization for area-preserving maps} In~\cite{GJ2} Gaidashev and Johnson construct simply connected domain $\DD_s
\subset \C^2$ and $\DD \subset \C$, and adapt the renormalization scheme
from~\cite{EKW2}. This renormalization scheme is defined on a neighborhood $\BB$
of $s_*\in \AAA(\DD_s)$, where $s_*$ corresponds to the Eckmann-Koch-Wittwer
fixed point.
There are analytic functions
\[
  F\mapsto \lambda_F\in (-\infty,0)
\]
and
\[
  F\mapsto \mu_F\in (0, \infty),
\]
called the rescaling, which are used to normalize $F$ (or $s$). The
renormalization operator $R$ is defined by
\[
  RF={\Lambda_F}^{-1}\circ F\circ F\circ \Lambda_F, \ {\rm where} \
  \Lambda_F:(x,u)\mapsto (\lambda_F x, \mu_F u).
\]

At the level of the generating functions, the renormalization operator $\RR: \BB
\subset \AAA(\DD_s) \mapsto \AAA(\DD_s)$ is defined as follows (see~\cite{GJ2}):
\begin{eqnarray} \label{Rs}
\RR[s](x,y)=\mu^{-1} s(z(x,y),\lambda y), 
\end{eqnarray}
where $z$ is the unique symmetric, that is, $z(x,y)=z(y,x)$, solution of 
\begin{equation}\label{midpoint}
s(\lambda x, z(x,y))+s(\lambda y, z(x,y))=0,
\end{equation}
and $\lambda$ and $\mu$ are fixed by the normalization conditions
$\RR[s](1,0)=0$, which implies $z(1,0)=z(0,1)=1$, and $\partial_1 \RR[s](1,0)=1$
which implies $\mu=\partial_1 z (1,0)$.

The results from~\cite{GJ2} which will be used in the sequel are
  collected in the following Theorem. All proofs in~\cite{GJ2} are done for the
operator $\RR$ in the neighborhood $\BB$ in the space $\AAA(\DD_s)$, however, it
is proved in~\cite{GJ2} that the map
\[
  \II:s \mapsto (y_s,s \circ h_s),
\]
where $y_{s}$ is the unique solution of $0=u+s(y,x)$, and
$h_{s}(x,u)=(x,y_s(x,u))$, is an analytic embedding, and the set $\FF:=\II(\BB)
\subset {\rm Cons}(\DD)$ is a Banach submanifold of the space $O(\DD)$ of
functions $F: \DD \mapsto \C^2$ holomorphic on $\DD$, continuous on $\partial
\DD$. The following is a Theorem from~\cite{GJ2}, reformulated for the
submanifold $\FF$.

\comm{ A map $F\in \text{Cons}(D)$ is called \textit{renormalizable} if there
  exists a unique periodic point $(p_1,p_2)$ of period $2$ and a rectangle
\begin{equation}\label{shape}
B_0=[p_1-a_1,p_1+a_1]\times [p_2-a_2,p_2+a_2]\subset {[-1,1]}^{2}
\end{equation}
such that
\begin{equation}\label{P1}
F^2(B_0)\cap B_0\ne \emptyset,
\end{equation}
and 
\begin{equation}\label{P2}
F(B_0)\cap B_0=\emptyset.
\end{equation}
The collection of renormalizable maps is denoted by $\text{Cons}_0(D)$. The
renormalization operator $R:\text{Cons}_0(D)\to \text{Cons}(D)$ is defined as
follows. Given a map $F\in \text{Cons}_0(D)$ let $B_0$ be the smallest rectangle
with the properties (\ref{shape}), (\ref{P1}) and (\ref{P2}). Then
\[
  RF={\Lambda_F}^{-1}\circ F\circ F\circ \Lambda_F,
\]
where $\Lambda_F: \mathbb{R}^2\to \mathbb{R}^2$ is the affine map
\[
  \Lambda_F:(x,u)\mapsto (-\lambda_F x+h_1, \mu_F u+h_2)
\]
with $\Lambda_F(D)=B_0$ and $\lambda_F, \mu_F>0$. A map $F\in \text{Cons}_0(D)$
is infinitely renormalizable if $R^{n}F$ is defined for all $n\ge 1$.

The following Theorem is proved in~\cite{GJ2}. See also Koch etc.
}

\begin{theorem}\label{hyp}
  There exists $F_*\in \FF$, which is generated by a symmetric function $s^*$: $s_1^*(x,y)=s_1^*(y,x)$; such that
  \begin{itemize}
  \item[1)] $F_*$ is a hyperbolic fixed point of the renormalization operator $R:   \FF \mapsto {\rm Cons(\DD)}$.
  \item[2)] $F_*$ has a two-dimensional unstable manifold in $\FF$.
  \item[3)] $F_*$ has a codimension two stable manifold $W^s(F_*)$ in $\FF$.
  \item[4)] $F^*$ has a codimension three strong stable manifold 
    $W^{ss}(F_*)\subset W^s(F_*)$.
  \item[6)] There exist a distance function $d$ on $W^{ss}(F_*)$, and $\nu<0.126$
    such that for every $F, \tilde{F}\in W^{ss}(F_*)$
    \[
      d(RF,R\tilde{F})\le \nu \cdot d(F,\tilde{F}).
    \]
    In particular,
    \[
      d(R^{n}F,F_*)\le \nu^{n}\cdot d(F,F_*).
    \]
  \item[7)] The one dimensional family defined by $F_t=\phi_t^{-1}\circ F_*\circ
    \phi_t$, where $\phi_t:\DD \to \phi_t(\DD) \subset \mathbb{R}^2$ is the
    diffeomorphism defined by
    \[
      \phi_t(x,u)=\left(x+tx^2, \frac{u}{1+2tx}\right),
    \]
    for $|t|$ sufficiently small, is contained in the stable manifold $W^s(F_*)$
    and is transversal to the strong stable manifold $W^{ss}(F_*)\subset W^s(F_*)$ and intersects it only in $F_*$.
  \item[8)] The map $\Lambda_F$ depends analytically on $F$.    
  \end{itemize}
\end{theorem}
\begin{remark}
We would like to emphasize that parts $1)$, $2)$ and $3)$, as well as the bound on the contraction rate $\nu$ in this Theorem have been proved using computer assistance.
\end{remark}

For the renormalization fixed point $F_{*}$ we will use the simpler notation
$\Lambda_{*}\coloneqq \Lambda_{F_{*}}$ and
$\lambda_{*}\coloneqq \lambda_{F_{*}}$, $\mu_{*}\coloneqq\mu_{F_{*}}$.

\subsection{The ratchet phenomenon}\label{subsectionratchet}
Twist maps have a property called the \textit{ratchet phenomenon}. It means that
for any twist map \( F(x,u) = (y, v) \) satisfying
\[
  \frac{\partial y}{\partial u} \ge a > 0
\]
there are horizontal cones \( \Theta_{h} \) and vertical cones
\( \Theta_{v} \) such that if \( p^{\prime}\in p+\Theta_{v} \) then $
F(p^{\prime})\in F(p)+\Theta_{h} $ and that the angle of the cones depend only
on $ a $, see e.g. Lemma 12.1 of \cite{Go}. The same is true for negative twist
maps with
\[
  \frac{\partial y}{\partial u} \le a < 0 \, .
\]
More precisely a positive twist map maps any point \(p^{\prime}\) in the half
cone \(p + \Theta_v^{\pm}\) into the half cone \(F(p) + \Theta_h^{\pm}\) and a
negative twist map maps any point \(p^{\prime}\) in \(p+\Theta_v^{\pm}\) into
\(F(p) + \Theta_h^{\mp}\). 

\begin{figure}[H]
  \centering
  \begin{tikzpicture}
    \node at (0, 0) {\includegraphics{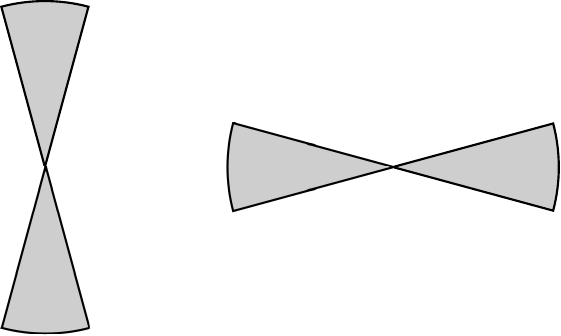}};
    \node at (-3.98, 0) {\(\bullet\)};
    \node at (1.9, 0) {\(\bullet\)};
    \node at (-3.8, 1.5) {\(\bullet\)};
    \node at (3.3, 0.15) {\(\bullet\)};
    \node at (-4.5, 0) {\(p\)};
    \node at (1.9, -0.6) {\(F(p)\)};
    \node at (-3.2, 1.5) {\(p^{\prime}\)};
    \node at (3.3, 0.8) {\(F(p^{\prime})\)};
    \node at (-4, 3.2) {\(\Theta_v^+\)};
    \node at (-4, -3.2) {\(\Theta_v^-\)};
    \node at (-1.3, 0) {\(\Theta_h^-\)};
    \node at (5.2, 0) {\(\Theta_h^+\)};
    \draw[->, >=stealth] (-3.4, 0) to node[above] {\(F\)} (-2, 0);
  \end{tikzpicture}
  \caption{The ratchet phenomenon for a positive twist map. A negative twist map
    reverses the signs.}
  \label{fig:ratchet}
\end{figure}
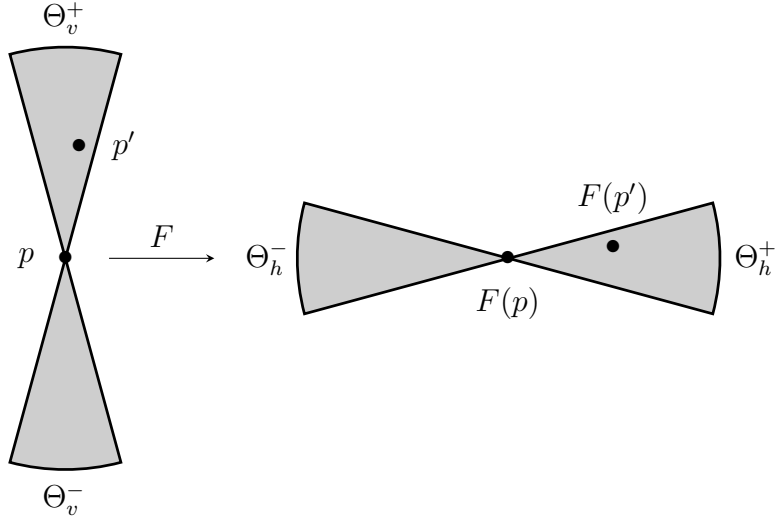

\section{The Cantor set}
\label{cantor}

In this section we recall the construction of the invariant Cantor set for
infinitely renormalizable maps. As in dimension one, it is a Cantor set on which the map acts like the dyadic adding machine. The construction is done via an
iterated function system.

We will use two kinds of notation
$$
\psi_0^F=\Lambda_F: D \to \mathbb{R}^2, \quad
\psi_1^F=F\circ\Lambda_{F}:D\to \mathbb{R}^2,
$$
as well as
$$
\psi_0^n=\Lambda_{R^{n-1}F}: D \to \mathbb{R}^2, \quad
\psi_1^n=R^{n-1}F\circ\Lambda_{R^{n-1}F}:D\to \mathbb{R}^2.
$$
Here $D$ denotes the real slice of $\DD$. Observe,
$$
R^n F=(\psi_0^n)^{-1}\circ R^{n-1}F \circ  R^{n-1}F \circ \psi_0^n.
$$
The convergence of $R^nF\to F_*$ and Theorem \ref{hyp}(6) immediately implies

\begin{lemma}
  \label{limitpsi}
  For every $F \in W^s(F_*)$ and any $k \ge 0$,
  $$
  \lim_{n\to\infty} ||\psi_{i}^{R^n F}-\psi_{i}^{F_*}||_{C^k}=0.
  $$
  for $i = 0,1$.
\end{lemma}

The above Lemma shows a crucial difference between conservative and dissipative
infinitely renormalizable maps: in the conservative case, the rescaling maps
converge to a diffeomorphism. In the dissipative case, the corresponding
rescaling $\psi_1^n$ converge to a degenerate map, a map with one-dimensional
image.

The construction of the Cantor set in the conservative case is exactly the same
as in the dissipative case. The difference lies in the asymptotic behavior of
the rescalings.

Let
$$
\Psi^F_{00}= \psi^1_0\circ \psi^2_0, \quad \Psi^F_{01}= \psi^1_0\circ \psi^2_1,
\quad \Psi^F_{10}=\psi^1_1\circ\psi^2_0,\quad \dots.
$$
and, proceeding this way, construct, for any $w=(w_1, \dots, w_n)\in\{0,1\}^n$,
$n\ge 1$, the maps
$$
 \Psi^F_w = \psi^1_{w_1}\circ\dots\circ \psi^n_{w_n}:\Dom(R^nF)\to \Dom(F).
$$
The transformations $\Psi^F_w$ will be referred to as the {\it  renormalization microscope}.

\begin{lemma}[Lemma 2.2 of \cite{GJM}]
  \label{bob1}
  For $F\in \FF$ there are analytically defined simply connected domains
  $B_0(F)\subset D$ and $B_1(F)=F(B_0(F)) \subset D$ such that
  \begin{equation}
    \label{P2}
    B_1(F)\cap B_0(F)=\emptyset,
  \end{equation}
  and
  \begin{equation}
    \label{P1}
    F^2(B_0(F))\cap B_0(F)\ne \emptyset.
  \end{equation}
  Moreover,
  $$
  \psi_i(B_0(RF)\cup B_1(RF))\subset B_i(F),
  $$
  with $i\in \{0,1\}$.

  The sets $B_0(F)$ and $B_1(F)$ contain a period $2$ hyperbolic
  orbit of $F$, $p^F_0 \in B_0(F)$, $p^F_1 = F(p^F_0) \in B_1(F)$, which is the
  unique such orbit in $B_0(F) \cup B_1(F)$.

  Additionally, the set $B_0$ can be chosen so that $T(B_0)=B_0$,
  $\Lambda_F(F(B_1)) \subset B_0$.
\end{lemma}
\begin{remark}\label{B0ellipse}
  There are likely many choices of the initial piece $B_0(F)$, however, we will always assume that this piece is an ellipse. A computer-assisted construction  $B_0(F)$  as an ellipse has been provided in Lemma $7.1$ in ~\cite{GJ3}. Additionally, it is shown in ~\cite{GJ3} that the pieces $B_0(F)$ and $B_1(F)$  intersect a horizontal line.
\end{remark}

\begin{proposition}[Proposition 2.3 of \cite{GJM}]
  \label{esti}
  There exists a neighborhood $\UU$ of $F_*$ and $0<\theta_1<\theta_2<1$ such
  that for $F\in \Wl$,
  \begin{equation}\label{Wl}
    \Wl=W^s(F_*) \cap \UU
  \end{equation}
  we have
  $$
  \theta_1^4 \cdot |v|\le |D\Psi^F_w(x,u) v|\le \theta_2^4 \cdot |v|
  $$
  for every $w\in \{0,1\}^4$ and $(x,u)\in B_0(F)\cup B_1(F)$. Moreover,
  \begin{equation}
    \label{crucialestimate}
    \frac{\theta_2\nu }{\theta_1}<1.
  \end{equation}
\end{proposition}

\begin{remark}
  The following estimates are obtained in \cite{GJ2}.
  \begin{equation}\label{theta1}
    \theta_1\ge 0.061
  \end{equation}
  \begin{equation}\label{theta2}
    \theta_2\le 0.249
  \end{equation}
  \begin{equation}\label{nu}
    \nu< 0.126
  \end{equation}
  \begin{equation}\label{ratio}
    \frac{\theta_2\nu }{\theta_1}< 0.515
  \end{equation}
  Note that these estimates also produce estimates of the values $\left|
    \lambda_{*} \right|$ and $\left| \mu_{*} \right|$. In fact according to
  Theorem 1 of \cite{GJ2} we have
  \begin{equation}\label{lambdamustar}
    \lambda_{*} \in [-0.24887681, -0.24887376], \quad
    \mu_{*} \in [0.061107811, 0.061112465].
  \end{equation}

We emphasize that all these bounds are computer-assisted.
\end{remark}

\begin{remark}
  \label{optcond}
  The estimate (\ref{crucialestimate}) plays a crucial role in the proof of the   Rigidity Theorem $\ref{rigidity}$. The numbers $\theta_2$, $\theta_1$ and
  $\nu$ are bounds on the maximal expansion and contraction in the
  renormalization microscope, and a bound on the spectral radius of
  renormalization. As it turns out, the Rigidity Theorem $\ref{rigidity}$ can be
  proved under the condition

  \vspace{1mm}
  
  $$
  \frac{{\rm maximal \ expansion} \cdot {\rm spectral \ radius} }{ {\rm maximal \
      contraction} }<1.
  $$

  \vspace{1mm}
  
  The derivative of renormalization at the fixed point is a compact operator. In
  particular, rigidity can be proved on a finite codimension subspace where the
  contraction is strong enough. The numerical estimates from \cite{GJ2} show
  that only the weakest stable direction is not strong enough. Luckily, this
  weakest stable direction corresponds to a one-dimensional family of
  analytically conjugated maps.
\end{remark}

The following Lemma follows directly from Proposition \ref{esti}.

\begin{lemma}
  \label{contracting}
  For every $F\in \Wl$ there exists $C>0$ such that for any word $w\in
  \{0,1\}^n$, $n\ge 1$,
  $$
  \| D\Psi^F_w\|\leq C \theta_2^n
  $$
  where $\theta_2<1$ is given in Proposition \ref{esti} and (\ref{theta2}).
\end{lemma}

Define the {\em pieces} of the $n^{\text{th}}$-{\em level} or
$n^{\text{th}}$-{\em scale} as follows. For each $w\in \{0,1\}^n$ let
$$
B^n_{w0}\equiv B^n_{w0}(F) = \Psi^F_w (B_{0}(R^nF))
$$
and
$$
B^n_{w1}\equiv B^n_{w1}(F) = \Psi^F_w (B_{1}(R^nF)).
$$

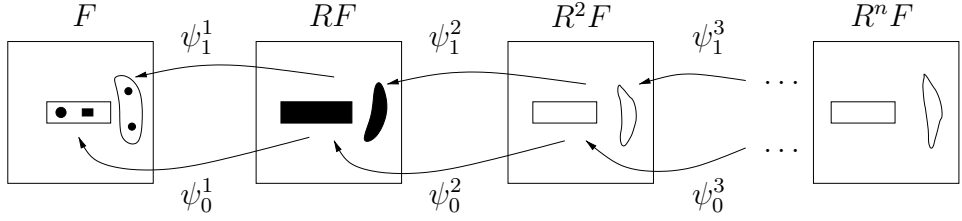
\begin{figure}[H]
\begin{center}
  \begin{tikzpicture}
    \node at (0, 0) {\includegraphics[width=\textwidth]{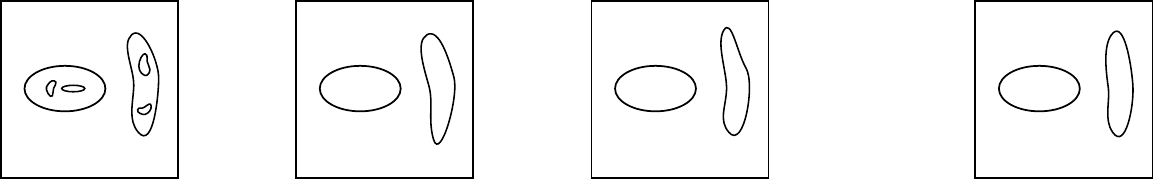}};
    \node at (4, -0.5) {$\dots$};
    \node at (4, 0.5) {$\dots$};
    \node at (6.3, 1.4) {$R^nF$};
    \node at (1.35, 1.4) {$R^2F$};
    \node at (-2.5, 1.4) {$RF$};
    \node at (-6.3, 1.4) {$F$};

    \draw[->, >=stealth] (3.8, 0.7) to[out=155, in = 25] node[shift={(0,0.4)}] {$\psi_{1}^{3}$} (2.2, 0.8);
    \draw[->, >=stealth] (3.8, -0.7) to[out=205, in = -45] node[shift={(0,-0.4)}] {$\psi_{0}^{3}$} (1.3, -0.5);
    \draw[->, >=stealth] (1, 0.7) to[out=155, in = 25] node[shift={(0,0.4)}] {$\psi_{1}^{2}$} (-1.7, 0.8);
    \draw[->, >=stealth] (1, -0.7) to[out=205, in = -45] node[shift={(0,-0.4)}] {$\psi_{0}^{2}$} (-2.6, -0.5);
    \draw[->, >=stealth] (-2.9, 0.7) to[out=155, in = 25] node[shift={(0,0.4)}] {$\psi_{1}^{1}$} (-5.5, 0.8);
    \draw[->, >=stealth] (-2.9, -0.7) to[out=205, in = -45] node[shift={(0,-0.4)}] {$\psi_{0}^{1}$} (-6.4, -0.5);
  \end{tikzpicture}
\caption{The Renormalization Microscope.}\label{micro}
\end{center}
\end{figure}

The set of words $\{0,1\}^n$ can be viewed as the additive group of residues
$\mod 2^n$ by letting
$$
w \mapsto\sum_{k=0}^{n-1} w_{k+1} 2^k.
$$
Let $p\colon \{0,1\}^n\ra \{0,1\}^n$ be the operation of adding 1 in this group.
The following has been proved in \cite{GJM}.

\begin{lemma}
  \label{permute}
  For every $F\in W^s(F_*)$ and $n\ge 1$
  \begin{enumerate}
  \item[1)] The above families of pieces are nested:
    $$
    B^n_{w\nu}\subset B^{n-1}_w, \quad w\in \{0,1\}^{n}, \ \nu\in \{0,1\}.
    $$
  \item[2)] The pieces $B^n_w$, $w\in \{0,1\}^{n+1}\setminus \{1^{n+1}\}$, are
    pairwise disjoint.
  \item[3)] Under $F$, the pieces are permuted as follows.
    $$
    F(B^n_w) = B^n_{p(w)},
    $$ 
    unless $p(w) = 0^{n+1}$. If $p(w)=0^{n+1}$, then $F(B^n_w) \cap
    B^n_{0^{n+1}}\ne \emptyset$.
  \end{enumerate}
\end{lemma}

The union of all pieces of level $n$ will be denoted by $\BB^n$ or by $\BB^n_F$
whenever the dependence on $F$ has to be emphasize:
$$\BB^n \equiv \BB^n_F= \bigcup_{w \in \{0,1\}^{n+1}} B^n_w.$$
 
Lemma \ref{contracting} implies:

\begin{lemma}
  \label{boxes shrink}
  For every $F\in \Wl$ there exists $C>0$ such that for all $w\in
  \{0,1\}^{n+1}$, $ \diam B^n_w\leq C \theta_2^n$.
\end{lemma}

Let
$$
\CC_F= \bigcap_{n=1}^\infty \bigcup_{w\in \{0,1\}^n} B^n_w.
$$
Let us also consider the {\it dyadic group} $ \{0,1\}^\N = \lim_{\leftarrow}
\{0,1\}^n$. For simplicity we will use the notation $\Omega \coloneqq
  \{0,1\}^{\N}$. The elements of $\Omega$ are infinite sequences $(w_1
w_2\dots)$ of symbols $0$ and $1$ that can be also represented as formal power
series
$$
w \mapsto\sum_{k=0}^\infty w_{k+1} 2^k.
$$
The natural numbers $\N$ are
embedded into $\Omega$ as finite series. The {\it adding machine} $p:
\Omega\ra \Omega$ is the operation of adding $1$ in this group. The
discussion above yields:

\begin{theorem}
  \label{adding machine}
  For every $F\in \Wl$ the set $\CC_F$ is an invariant Cantor set. The map $F$
  acts on $\CC_F$ as the dyadic adding machine $p: \Omega\ra
  \Omega$. The conjugacy between $p$ and $F|_{\CC_F}$ is given by the
  following homeomorphism $h_F: \Omega\ra \CC_F$:
  \begin{equation}
    \label{conjh}
    h_F: w = (w_1w_2\dots) \mapsto \bigcap_{n=1}^\infty B^n_{w_1\dots w_{n+1}}. 
  \end{equation}
  Furthermore, $\CC_F$ has Lebesgue measure zero with
  $$
  \text{HD}(\CC_F)\le -\frac{\log 2}{\log \theta_2}\le 0.795.
  $$
\end{theorem}

The invariant Cantor sets $\CC_F$ are the counterpart of the period doubling
Cantor sets in one-dimensional dynamics and strongly dissipative H\'enon-like
maps, see \cite{CLM, GST, Mi}. The dynamics of $F$ restricted to $\CC_F$ is
conjugated to the adding machine. The adding machine is uniquely ergodic. Let
$\mu$ be the unique invariant measure of $F$ restricted to $\CC_F$:
\begin{equation}
  \label{eq:cantorSetMeasure}
  \mu(B^n_{w})=\frac{1}{2^{n+1}}.
\end{equation}

The proof of the following theorem appears in \cite{GJM}.

\begin{theorem}
  \label{expo}
  The measure $\mu_F$ of every map $F\in W^s(F_*)$ has a single characteristic
  exponent, $\chi=0$.
\end{theorem}

The most important result concerning the period doubling Cantor sets for
area-preserving maps is their rigidity in the local stable manifold
  $\Wl$, see Theorem 4.1 of \cite{GJM}:

\begin{theorem}({\it \underline{Rigidity}}.)
  \label{rigidity}
  Let $F$, $G$ be any two maps in $\Wl$. Then there exists
  $\alpha>0.237$, such that the topological conjugacy between $F
  \arrowvert_{\CC_F}$ and $G \arrowvert_{\CC_G}$ extends to a $C^{1+\alpha}$
  diffeomorphism of neighborhoods of the Cantor sets.
\end{theorem}

\section{Expansion and contraction in the pieces}
\label{Lyapunov}

\subsection{Vanishing hyperbolicity of periodic orbits}
Let $D \subset \R^2$ be the domain of analyticity of $F$, and suppose that $p$
is a point whose orbit $\mathcal{O}(p) \subset D$. Recall, the definition of the
upper Lyapunov exponent of $(p,v) \in D \times \R^2$ with respect to $F$:
$$
\chi(p,v;F) \equiv \overline{\lim}_{i \rightarrow \infty} {1 \over i}
\ln\left[ \Arrowvert DF^i(p) v \Arrowvert\right],
$$
where $\Arrowvert \cdot \Arrowvert$ is some norm in $\R^2$. The maximal and
minimal Lyapunov exponent of $p \in D$ with respect to $F$ are defined as
\begin{align*}
  \chi_{+}(p;F) \equiv \sup_{||v||=1} \chi(p,v;F), \\
  \chi_{-}(p;F) \equiv \inf_{||v||=1} \chi(p,v;F) \\
\end{align*}
respectively.

The following Lemma about the existence of hyperbolic fixed points for maps in a
small neighborhood of the renormalization fixed point map $F^*$ is a restatement
of a result from \cite{GJ1}. The proof of the Lemma is computer-assisted (see
\cite{GJ1}).

\begin{lemma}
  \label{hyp_point}
  Let $\Wl$ be as in \eqref{Wl}. Every map $F \in \Wl$ possesses a hyperbolic
  fixed point $p^F \in D$, such that
  \begin{itemize}
  \item[1)] $\pi_x p^F \approx 0.5776$, and $\pi_u p^F=0$, where $\pi_{x,u}$ are
    projections on the $x$ and $u$ coordinates;
  \item[2)] $D F(p^F)$ has two negative eigenvalues.
    \[
      e_-^F \approx 0.486, \quad e_+^{F} = \frac{1}{e_-^{F}}
    \]
    corresponding to the following two eigenvectors:
    \[
      {s}^F \approx [1.0,0.78], \quad {\rm and} \quad {u}^F=T({s}^F).
    \]
  \end{itemize}
\end{lemma}

\bigskip

An immediate consequence of this Lemma is existence of hyperbolic $2^n$-th
periodic orbits for maps in $\Wl$. Let $\mathcal{O}_{n}(F)$ denote such $2^n$-th
periodic orbit for $F \in \Wl$, specifically:
$$
\mathcal{O}_{n}(F)=\cup_{i=0}^{2^n-1} F^{i}(\Psi^{F}_{0^n}(p^{F_n})),
$$
where $p^{F_n}$ is the fixed point of
$$F_n \equiv R^n[F] \in \Wl.$$
We will also denote
$$
p_{0^n}^F=\Psi^{F}_{0^n}(p^{F_n}),\quad p_{w}^F= F^{\sum_{i=1}^n w_i 2^{i-1}}(p_{0^n}^F).
$$

Consider the stable and unstable invariant direction fields on
$\mathcal{O}_{n}(F)$. At every point $p_w^F$, $w \in \{0,1\}^n$ of
$\mathcal{O}_{n}(F)$, these directions are given by
\begin{eqnarray}
  \label{vec_u} {u}_w^F&=&D \left( F^{\sum_{i=1}^n w_i 2^{i-1}} \circ \Psi^{F}_{0^n}\right)(p^{F_n}) {u}^{F_n}= D \Psi^F_w(p^{F_n}) {u}^{F_n},\\ 
  \label{vec_s} {s}_w^F&=&D \left( F^{\sum_{i=1}^n w_i 2^{i-1}} \circ \Psi^{F}_{0^n}\right)(p^{F_n}) {s}^{F_n}= D \Psi^F_w(p^{F_n}) {s}^{F_n}. 
\end{eqnarray}
The subspaces defined by the span of these vectors will be denoted by

$E_{-}(p^{F}_{w})\coloneqq \mathrm{span}\{s^{F}_{w}\}$ and
$E_{+}(p^{F}_{w})\coloneqq \mathrm{span}\{u^{F}_{w}\}$ respectively.


\begin{lemma}
  \label{orbit_bounds}
  The set $\cup_{n=0}^\infty \mathcal{O}_{n}(F) \cup \CC_F$ is in the set of regular points
  for $F$, specifically,
  The decomposition
    $$
    \R^2=E_-(p_w^F) \bigoplus E_+(p_w^F), \ w \in \{0,1\}^n 
    $$
    is invariant under
    $$
    D F: D \times \R^2 \mapsto \R^2.
    $$

    The Lyapunov exponents
    $$
    \chi_-(p;F)\equiv -\chi_+(p;F)= {\ln{|e_-^{F_n}|} \over 2^n}, \quad p \in
    \mathcal{O}_{n}(F),
    $$
    where $e_-^F$ is as in Lemma $\ref{hyp_point}$, exist, are $F$-invariant,
    and
    $$
    \lim_{i \rightarrow \infty} {1 \over i} \ln\left\{  { \Arrowvert D F^i(x) v
        \Arrowvert  \over \|v\| }\right\}= \chi_\pm(p;F),
    $$
    for all $v \in E_\pm(p^F_w) \setminus \{0\}$;

\end{lemma}
\begin{proof}
  Let $i=q 2^n +k$, $k=2^{j_1}+2^{j_2}+ \ldots + 2^{j_m} < 2^n$, then
  \begin{eqnarray*}
    DF^i(p_{0^n}^F) {s}_{0^n}^F&=&DF^{k+q 2^n}(p_{0^n}^F) {s}_{0^n}^F=DF^k(F^{q 2^n}(p_{0^n}^F))  \cdot DF^{q 2^n}(p_{0^n}^F) {s}_{0^n}^F\\
    &=&DF^k(p_{0^n}^F) \cdot D \left( \Psi^{F}_{0^n} \circ  F^q_n \circ \left(\Psi^{F}_{0^n}\right)^{-1} \right)(p_{0^n}^F) {s}_{0^n}^F   \\
    &=&DF^k(p_{0^n}^F) \cdot D \Psi^{F}_{0^n} \left(  F^q_n \circ \left( \Psi^{F}_{0^n}\right)^{-1} (p_{0^n}^F) \right) \\
    &\phantom{=}&\phantom{DF^k(p_{0^n})}  \cdot DF^q_n \left( \left(\Psi^{F}_{0^n}\right)^{-1}\!\!(p_{0^n}^F)  \right)   \cdot \! D \left(\Psi^{F}_{0^n}\right)^{-1}\!\! (p_{0^n}^F)  {s}_{0^n}^F \\
    &=&DF^k(p_{0^n}^F) \cdot \! D \Psi^{F}_{0^n} (p^{F_n}) \cdot DF^q_n (p^{F_n})  \cdot \! D \left(\Psi^{F}_{0^n}\right)^{-1}\!\!(p_{0^n}^F)  {s}_{0^n}^F \\
    &=&DF^k(p_{0^n}^F) \cdot D \Psi^{F}_{0^n} (p^{F_n}) \cdot DF^q_n (p^{F_n}) {s}^{F_n}\\
    &=&DF^k(p_{0^n}^F) \cdot D \Psi^{F}_{0^n} (p^{F_n}) \cdot (e_-^{F_n})^{q} {s}^{F_n}\\
    &=&  (e_-^{F_n})^{q} DF^k(p_{0^n}^F) \cdot  {s}_{0^n}^F.
  \end{eqnarray*}

  {Let $C_{n}$ and
    $c_{n}$ denote} upper and lower bounds on the derivative norm of $F$ on
  $\mathcal{O}_{n}(F)$ {respectively}. Then
  \begin{equation}
    c_n^k |e_-^{F_n}|^q \Arrowvert{s}_{0^n}^F\Arrowvert  \le \Arrowvert DF^i(p_{0^n}^F) {s}_{0^n}^F \Arrowvert  \le C_n^k |e_-^{F_n}|^q \Arrowvert{s}_{0^n}^F\Arrowvert ,
  \end{equation}
  and,
  \begin{eqnarray*}
    \overline{\lim}_{i \rightarrow\infty} {1 \over i} \ln \left[ {\Arrowvert DF^i(p_{0^n}^F) {s}_{0^n}^F \Arrowvert \over \| {s}_{0^n}^F \|} \right] &\le& \lim_{i \rightarrow \infty} {k \over i} \ln C_n+ {q \over i} \ln\{|e_-^{F_n}|\} \\
    &=&{\ln\{|e_-^{F_n}|\} \over 2^n},\\
    \underline{\lim}_{i \rightarrow\infty} {1 \over i} \ln \left[ {\Arrowvert DF^i(p_{0^n}^F) {s}_{0^n}^F \Arrowvert \over \| {s}_{0^n}^F \|} \right] &\ge& \lim_{i \rightarrow \infty} {k \over i} \ln c_n+ {q \over i} \ln\{|e_-^{F_n}|\} \\ 
    &=&{ \ln\{|e_-^{F_n}|\} \over 2^n},
  \end{eqnarray*}
  therefore, the limit
  \[
    \lim_{i \rightarrow\infty} {1 \over i} \ln \left[ {\Arrowvert
        DF^i(p_{0^n}^F) {s}_{0^n}^F \Arrowvert \over \| {s}_{0^n}^F \|} \right]
  \]
  exists, and is equal to
  \[
    \chi_-(p;F) \equiv { \ln\{|e_-^{F_n}|\} \over 2^n}.
  \]

  A similar computation demonstrates that
  \[
    \lim_{i \rightarrow\infty} {1 \over i} \ln \left[ {\Arrowvert
        DF^i(p_{0^n}^F) {u}_{0^n}^F \Arrowvert \over \| {u}_{0^n}^F \|}
    \right]={ \ln\{|e_+^{F_n}|\} \over 2^n}=-\chi_-(p;F).
  \]

  \bigskip


\end{proof}

\subsection{Simple bounds on expansion and contraction in the pieces}
In this subsection we will obtain some initial bounds on expansion and
contraction in pieces of a fixed level. These bounds will be used to prove
positive measure of {bounded and} unbounded geometry.

\begin{lemma}\label{expansion_bounds}
  There is a constant $C>0$ such that for any $D>0$ satisfying
  \begin{equation}
    \label{Dcond1}  D  \ge  8 {n \over n+3} (\ln |\lambda_*| -\ln \mu_*)+ {8 \over n+3} \ln \| D F_* \arrowvert_{B_{w}} \|, \\
  \end{equation}
  for all $n \in \N$ and $w \in \{0,1\}$, the following bound holds
  \begin{equation}\label{bound_upper}
    \| D F^i_* \arrowvert_{B_{0^{n-1}w}^{n-1}} \| \le  C   e^{D {n \over 2^n} i}
  \end{equation}
  for all $n \in \N$ and  $ 0 \le i < 2^{n-1}$
\end{lemma}

\begin{proof}  
  Clearly,  there exists $C>0$, such that
  \[
    \| D F^{k}_* \arrowvert_{B_{0w}^1} \| \le  C, \quad w \in \{0,1\}
  \]
  for all $0 \le k <2$. 
  Next, let $N \in \N$, $N>1$ and assume that the claim is true for all $1<n \le N$, and prove it for $n=N+1$.

  Let $i=2^{N-1}+k \ge 2^{N-1}$, $k\le 2^{N-2}$, that is, first we allow for the
  orbits of length up to $2^{N-1}+2^{N-2}$ only, and consider $D F^i
  \arrowvert_{B_{0^{N}w}^{N}}$.
  \begin{eqnarray*}
    \| D F_*^{i} \arrowvert_{B_{0^{N}w}^{N}} \|&=& \| D \left( F_*^{k} \circ F_*^{2^{N-1}} \right) \arrowvert_{B_{0^{N}w}^{N}} \|  \\
    &=& \| D \left( F^k_* \circ \Lambda_{*}^{N-1} \circ  R^{N-1} F_*  \circ   \Lambda_{*}^{-N+1}  \right) \arrowvert_{B_{0^{N}w}^{N}} \| \\
    & \le & {|\lambda_*|^{N-1} \over \mu_*^{N-1} } \| D F_* \arrowvert_{B_{0w}^1} \| \| D F^{k}_* \arrowvert_{B^{N}_{0^{N-1}1w}} \|.
  \end{eqnarray*}
  We remark at this point, that $k \le 2^{N-2} \le 2^{N-1}$ and that
  $B^N_{0^{N-1}1w} \subset B^{N-1}_{0^{N-1}1}$, and thus to estimate $\| D
  F^{k}_* \arrowvert_{B^{N}_{0^{N-1}1w}} \|$ we can use the (assumed) bound
  \eqref{bound_upper} with $n=N$:
  \begin{align}
    \nonumber  \| D F_*^{i} \arrowvert_{B_{0^{N}w}^{N}} \|& \le {|\lambda_*|^{N-1} \over \mu_*^{N-1} }  \| D F_* \arrowvert_{B_{0w}^1} \|  C e^{D {N \over 2^N}\left(i-2^{N-1}\right) } \\
    \nonumber &\le {|\lambda_{*}|^{N-1} \over \mu_*^{N-1}}   \| D F_* \arrowvert_{B_{0w}^1} \|  e^{-D {N \over 2}+D \left({N \over 2^N}-{N+1 \over 2^{N+1}}   \right) i} C e^{ D {N+1 \over 2^{N+1}} i } \\
    \nonumber &\le {|\lambda_{*}|^{N-1} \over \mu_*^{N-1}}  \| D F_* \arrowvert_{B_{0w}^1} \|  e^{-D {N \over 2}+D \left({N \over 2^N}-{N+1 \over 2^{N+1}}   \right) (2^{N-1}+2^{N-2})} C e^{ D {N+1 \over 2^{N+1}} i } \\
    \label{Dcondd1} &\le {|\lambda_{*}|^{N-1} \over \mu_*^{N-1}}  \| D F_* \arrowvert_{B_{0w}^1} \| e^{-D {N+3 \over 8}} C e^{ D {N+1 \over 2^{N+1}} i }.
  \end{align}
  
  Next, we consider $i=2^{N}-k$, $1\leq k \leq 2^{N-2}$, that is, we consider
  the orbits of length between $2^{N-1}+2^{N-2}$ and $2^{N} - 1$:
  \begin{eqnarray*}
    \| D F_*^{i} \arrowvert_{B_{0^{N}w}^{N}} \|&=& \| D \left( F_*^{-k} \circ F_*^{2^{N}} \right) \arrowvert_{B_{0^{N}w}^{N}} \|  \\
    &=& \| D \left( F^{-k}_* \circ \Lambda_{*}^{N} \circ  R^{N} F_*  \circ   \Lambda_{*}^{-N}  \right) \arrowvert_{B_{0^{N}w}^{N}} \| \\
    & \le & {|\lambda_*|^{N} \over \mu_*^{N} }  \| D (T \circ F^{k}_* \circ T) \arrowvert_{\Lambda_*^N(F_*(B_w^0))} \| \| D F_* \arrowvert_{B_{w}^0} \|  \\
    & \le & {|\lambda_*|^{N} \over \mu_*^{N} }  \| D F^{k}_* \arrowvert_{T(\Lambda_*^N(F_*(B_w^0)))} \| \| D F_* \arrowvert_{B_{w}^0} \|.
  \end{eqnarray*}

  Notice that $T(\Lambda_*^N(F_*(B_0^0)))=T(B_{0^N1}^N)$, and since
  $T(B_{0^N}^{N-1})=B_{0^N}^{N-1}$ by Lemma~\ref{bob1}, we get that
  $T(B_{0^N1}^N) \subset B_{0^N}^{N-1}$. Again, by Lemma~\ref{bob1},
  $\Lambda_*(F_*(B_1^0)) \subset B_0^0$, therefore $T(\Lambda_*^N(F_*(B_1^0)))
  \subset T(\Lambda_*^{N-1}(B_0^0))=B^{N-1}_{0^N}$. We conclude that for both
  $w=0$ and $w=1$,
  \[
    {T(\Lambda_*^N(F_*(B_w^0)))} \subset B^{N-1}_{0^N},
  \]
  and we can estimate $\| D F^{k}_* \arrowvert_{T(\Lambda_*^N(F_*(B_w^0)))} \|$
  using~\eqref{bound_upper} with $n=N$ once again:
  \begin{eqnarray}
    \nonumber  \| D F_*^{i} \arrowvert_{B_{0^{N}w}^{N}} \| & \le & {|\lambda_*|^{N} \over \mu_*^{N} }  C  \| D F_* \arrowvert_{B_{w}^0} \| e^{D {N \over 2^N} k} \\
    \nonumber   & \le & {|\lambda_*|^{N} \over \mu_*^{N} }    \| D F_* \arrowvert_{B_{w}^0} \| C e^{D {N \over 2^N} (2^N-i)} \\
    \nonumber  & \le & {|\lambda_*|^{N} \over \mu_*^{N} }    \| D F_* \arrowvert_{B_{w}^0} \| e^{D N -D \left( {N \over 2^N}+{N+1 \over 2^{N+1}}\right) \left(2^{N-1}+2^{N-2}  \right) } C e^{D {N+1 \over 2^{N+1}} i} \\
    \label{Dcondd2}  & \le & {|\lambda_*|^{N} \over \mu_*^{N} }   \| D F_* \arrowvert_{B_{w}^0} \|  C e^{-D {N+3 \over 8}} e^{D {N+1 \over 2^{N+1}} i}.
  \end{eqnarray}
  The conditions
  $$ {|\lambda_{*}|^{N-1} \over  \mu_*^{N-}} \| D F_*  \arrowvert_{B_{0w}^1} \| e^{-D {N+3 \over 8}} \le 1 \hspace{2mm} {\rm and} \hspace{2mm} {|\lambda_*|^{N}  \over \mu_*^{N} } \| D F_* \arrowvert_{B_{v}^0} \| e^{-D {N+3 \over 8}} \le 1$$
  are implied by $(\ref{Dcond1})$ and guarantee that (\ref{Dcondd1}) and
  (\ref{Dcondd2}) are bounded by $C e^{D {N+1 \over 2^{N+1}} i}$. The claim
  follows.
\end{proof}

Next we find a bound on contraction within the pieces from below.

\begin{lemma}\label{contraction_bounds}
  There is a constant $C>0$ such that for any $D>0$ satisfying
  \begin{equation}
    \label{eq:D_condition_lower}
    D > 8\frac{n}{n+3}\left( \ln|\lambda_{*}| - \ln\mu_{*} \right) - \frac{8}{n+3}\ln b
  \end{equation}
  where $b$ is a lower bound on
  \[
    \min_{w \in \{0,1\}, x \in B_w, v \in \R^2 \setminus \{0\}} {\| DF(x) v \| \over \|v\|},
  \]
  the following holds
  \begin{equation}
    \label{bound_lower}
    \| D F^{i}_* \arrowvert_{B_{0^{n-1}w}^{n-1}} v \| \ge Ce^{-D {n \over 2^n} i} \| v \|  
  \end{equation}
  for all $n \in \N$, $0 \le i <2^{n-1}$, $w\in\{0,1\}$ and all $v \in \R^{2}$.
\end{lemma}
  
\begin{proof}
  Again, there exists $C>0$, such that
  \[
    \| D F^{k}_* \arrowvert_{B_{0w}^{1}} v \| \ge C\|v\|
  \]
  for all $0 \le k <2$ and $v \in \R^{2}$.

  Let $N \in \N$, $N>1$, and assume that the claim is true for all $1<n \le N$.
  As before let $i=2^{N-1}+k \ge 2^{N-1}$, $k<2^{N-2}$, and consider the action
  of $D F^{i}_* \arrowvert_{B_{0^{N} w}^{N}}$ on $v \in \R^{2}$. Note that, like
  before, we have $B^{N}_{0^{N-1}1w}\subset B^{N-1}_{0^{N-1}1}$ and so we can
  use the (assumed) bound~(\ref{bound_lower}) with $n=N$ to bound the action of
  $DF^{k}_{*}\vert_{B^{N}_{0^{N-1}1w}}$.
  \begin{align*}
    \| D F^{i}_* \arrowvert_{B_{0^{N} w}^{N}} v \| & = \| D \left( F^{k}_* \circ F^{2^{N-1}}_*  \right) \arrowvert_{ B_{0^{N} w}^{N} } v \| \\
                                                   & = \| D \left( F^{k}_* \circ  \Lambda_{*}^{N-1} \circ R^{N-1} F_* \circ  \Lambda_{*}^{-N+1}  \right) \arrowvert_{B_{0^{N} w}^{N}} v \| \\
                                                   & \ge \left\{ \min_{(x,u) \in B_{0^{N-1}1w}^{N}, v \in \R^{2} \setminus \{0\}} { \| D F^{k}(x,u) v \| \over \|v\| }  \right\} {\mu_{*}^{N-1} \over |\lambda_*|^{N-1}} b  \| v\|\\
                                                   & \ge Ce^{-D {{N} \over 2^{N}}\left(i -2^{N-1}\right)} {\mu_{*}^{N-1} \over |\lambda_*|^{N-1}} b  \| v\|\\
                                                   & \ge Ce^{D\frac{N}{2}}e^{-D {{N} \over 2^{N}}i} {\mu_{*}^{N-1} \over |\lambda_*|^{N-1}} b  \| v\|\\
                                                   & \ge Ce^{D\frac{N}{2}}e^{-D \left( {{N} \over 2^{N}} - {{N+1} \over {2^{N+1}}} \right)i}e^{-D{{N+1} \over 2^{N+1}}i} {\mu_{*}^{N-1} \over |\lambda_*|^{N-1}} b  \| v\|\\
                                                   & \ge Ce^{D\frac{N}{2}}e^{-D \left( {{N} \over 2^{N}} -
                                                     {{N+1} \over {2^{N+1}}} \right)\left( 2^{N-1} + 2^{N-2} \right)}e^{-D{{N+1} \over 2^{N+1}}i} {\mu_{*}^{N-1} \over |\lambda_*|^{N-1}} b  \| v\|\\
                                                   & \ge e^{D{{N+3} \over 8}} {\mu_{*}^{N-1} \over |\lambda_*|^{N-1}} \hspace{0.4mm} b \hspace{0.4mm} e^{ -D{N+1 \over 2^{N+1}} i }  \| v\|,
  \end{align*}
  We therefore have that
  \[
    \| D F^{i} \arrowvert_{B_{0^{N} w}^{N}} v\| \geq Ce^{ - D {N+1 \over 2^{N+1}} i}\|v\|
  \]
  for all $0 \le i \leq 2^{N-1} + 2^{N-2}$ and $w = 0,1$ due
  to~\eqref{eq:D_condition_lower}.

  Now let $i = 2^{N} - k$ where we let $1\leq k \leq 2^{N-2}$. Following the
  proof of Lemma~\ref{expansion_bounds} we have
  \begin{align*}
    \| D F_*^{i} \arrowvert_{B_{0^{N}w}^{N}} v \| & = \| D \left(F_*^{-k} \circ F_*^{2^{N}} \right) \arrowvert_{B_{0^{N}w}^{N}} v\| \\
                                                  & = \| D \left( F^{-k}_* \circ \Lambda_{*}^{N} \circ  R^{N} F_* \circ \Lambda_{*}^{-N} \right) \arrowvert_{B_{0^{N}w}^{N}} v \| \\
                                                  & \geq {\mu_{*}^{N} \over |\lambda_{*}|^{N} } \left\{\min_{(x,u) \in \Lambda_*^N(F_*(B_w)), v \in \R^{2} \setminus \{0\}}{ {\| D (T \circ F^{k}_* \circ T) v \|} \over {\|v\|}} \right\} b \| v \|.
  \end{align*}
  Using the same argument as in the proof of Lemma~\ref{expansion_bounds} we can
  see that $T(\Lambda_{*}^{N}(F_{*}(B_{w})))\subset B^{N-1}_{0^{N}}$ allowing us
  to use~\eqref{bound_lower} to make the following estimate:
  \begin{align*}
    \| D F_*^{i} \arrowvert_{B_{0^{N}w}^{N}} v \| & \geq {\mu_{*}^{N} \over |\lambda_{*}|^{N} } \left\{\min_{(x,u) \in \Lambda_*^N(F_*(B_w)), v \in \R^{2} \setminus \{0\}}{ {\| D (T \circ F^{k}_* \circ T) v \|} \over {\|v\|}} \right\} b \| v \| \\
                                                  & \geq {\mu_{*}^{N} \over
                                                    |\lambda_{*}|^{N}}Ce^{-D\frac{N}{2^{N}}k}b\|v\| \\
                                                  & = {\mu_{*}^{N} \over
                                                    |\lambda_{*}|^{N}}Ce^{-D\frac{N}{2^{N}}(2^{N}-i)}b\|v\| \\
                                                  & = {\mu_{*}^{N} \over
                                                    |\lambda_{*}|^{N}}Ce^{-DN}e^{D\frac{N}{2^{N}}i}b\|v\| \\
                                                  & = {\mu_{*}^{N} \over |\lambda_{*}|^{N}}Ce^{-DN}e^{D\left(\frac{N}{2^{N}} + \frac{N+1}{2^{N+1}}\right)i}e^{-D\frac{N+1}{2^{N+1}}i}b\|v\| \\
                                                  & \geq {\mu_{*}^{N} \over
                                                    |\lambda_{*}|^{N}}Ce^{-DN}e^{D\left(\frac{N}{2^{N}}
                                                    + \frac{N+1}{2^{N+1}}\right)(2^{N-1}
                                                    + 2^{N-2})}e^{-D\frac{N+1}{2^{N+1}}i}b\|v\| \\
                                                  & = {\mu_{*}^{N} \over
                                                    |\lambda_{*}|^{N}}Ce^{D\frac{N+3}{8}}e^{-D\frac{N+1}{2^{N+1}}i}b\|v\|.
  \end{align*}
  Again due to~\eqref{eq:D_condition_lower} we see that we get
  \[
    \| D F_*^{i} \arrowvert_{B_{0^{N}w}^{N}} v \| \geq Ce^{-D {N+1 \over 2^{N+1}} i}\|v\|.
  \]
  The claim follows.
\end{proof}

To prove existence of both bounded and unbounded geometry on a set of positive
measure we require a good control of the second derivatives of
the maps.

The second derivative of $F_*$ is an element of the Banach space $BL(\R^2
\times \R^2, \R^2)$ of all continuous bilinear maps from $\R^2$ to $\R^2$. The
action of this bilinear form on a pair of vectors $(u,v)$ will be denoted by $u
\cdot D^2 F(x) \cdot v$.

\comment{
\textcolor{magenta}{The notation $L^{2}$ is typically used for square-integrable
  maps so I think we should use some other notation. As far as I know there's
  not really a standard notation for this but the notation $T^{k}_{l}(M)$ is
  pretty common to use for the bundle of $(k,l)$ tensors over a smooth manifold
  $M$. Now, the space of sections of $T^{k}_{l}(M)$ really doesn't have any
  standard notation but I prefer the notation $\mathcal{T}^{k}_{l}(M)$. With
  this notation the Hessian of $F_{*}$ would be an element of
  $\mathcal{T}^{2}_{1}(\mathbb{R}^{2})$. What do you think of this notation?}
}

We will now obtain upper bounds on the second derivative of $F_*$ following
closely the ideas of Lemma \ref{expansion_bounds}.

\begin{lemma}\label{expansion_bounds_second}
  There is a constant $C>0$ such that for any $D>0$ satisfying
  \begin{equation}
    \label{DDcond1}  D  \ge  {8 n \over n+3} (\ln |\lambda_*| -2 \ln \mu_*)+ {8 \over n+3} \ln \left( |\lambda_*|^{n} \| D F_* \arrowvert_{B_{w}} \|^2 +  \| D^{2} F_* \arrowvert_{B_{w}} \|\right), \\
  \end{equation}
  for all $n \in \N$ and $w \in \{0,1\}$, the following bound holds
  \begin{equation}\label{bound_upper_Second} \| D^2 F^i_* \arrowvert_{B_{0^{n-1}w}^{n-1}} \| \le  C   e^{D {n \over 2^n} i}
  \end{equation}
  for all $n \in \N$ and  $ 0 \le i < 2^{n-1}$.
\end{lemma}

\begin{proof}  
  Clearly, there exist $B>0$ and $A>0$, such that
  \[
    \| D F^{k}_* \arrowvert_{B_{0w}^1} \| \le  B, \quad     \| D^2 F^{k}_* \arrowvert_{B_{0w}^1} \| \le  A,  \quad w \in \{0,1\}
  \]
  for all $0 \le k <2$. Next, let $N \in \N$, $N>1$. We assume that the claim is
  true for all $1<n \le N$, and prove it for $n=N+1$.

  Let $i=2^{N-1}+k \ge 2^{N-1}$, $k\le 2^{N-2}$, that is, first we allow for the
  orbits of length up to $2^{N-1}+2^{N-2}$ only, and consider $D^2 F^i
  \arrowvert_{B_{0^{N}w}^{N}}$. Denote, for brevity $G=\Lambda_{*}^{N-1} \circ
  F_* \circ \Lambda_{*}^{-N+1}$, then
  \begin{align*}
    \| D^2 F_*^{i} \arrowvert_{B_{0^{N}w}^{N}} \|&= \| D^2 \left( F_*^{k} \circ F_*^{2^{N-1}} \right) \arrowvert_{B_{0^{N}w}^{N}} \|  \\
    &= \| D^2 \left( F^k_* \circ \Lambda_{*}^{N-1} \circ F_*  \circ   \Lambda_{*}^{-N+1}  \right) \arrowvert_{B_{0^{N}w}^{N}} \| \\
    &= \left\|\left( DG  \cdot \left( D^2 F^k_* \right) \circ G \cdot DG \arrowvert_{B_{0^{N}w}^{N}} \hspace{-2mm} + (D F^k_* \circ G) \cdot D^2 G \right) \arrowvert_{B_{0^{N}w}^{N}} \right\| \\
    & \le  {|\lambda_*|^{2(N-1)} \over \mu_*^{2(N-1)} }  \| D F_* \arrowvert_{B_{0w}^1} \|^2 \| D^2 F^{k}_* \arrowvert_{B^{N}_{0^{N-1}1w}} \| +  \\
    & \hspace{10mm}+ \| D F^{k}_* \arrowvert_{B^{N}_{0^{N-1}1w}} \| \|\Lambda_*^{N-1} \| \| D^2 F_* \arrowvert_{B^{1}_{0w}}\| \|\Lambda_*^{-N+1}\|^2 \\
    & \le  {|\lambda_*|^{2(N-1)} \over \mu_*^{2(N-1)} }  B^2 \| D^2 F^{k}_* \arrowvert_{B^{N}_{0^{N-1}1w}} \| + A  {|\lambda_*|^{N-1} \over \mu_*^{2(N-1)} }  \| D F^{k}_* \arrowvert_{B^{N}_{0^{N-1}1w}} \|
  \end{align*}
  We remark at this point, that $k \le 2^{N-2} \le 2^{N-1}$ and that
  $B^N_{0^{N-1}1w} \subset B^{N-1}_{0^{N-1}1}$, and thus to estimate $\| D^2
  F^{k}_* \arrowvert_{B^{N}_{0^{N-1}1w}} \|$ we can use the (assumed) bound with
  $n=N$. Additionally, since $D$ in the hypothesis is larger than $D$ in Lemma
  \ref{expansion_bounds}, we can use the upper bound on the first derivative
  from that Lemma, with $D$ as in \eqref{DDcond1}.
  \begin{align*}
    \| D^2 F_*^{i} \arrowvert_{B_{0^{N}w}^{N}} \|& \le{|\lambda_*|^{2(N-1)} \over \mu_*^{2(N-1)} }  B^2  C e^{ D {N \over 2^N}\left(i-2^{N-1}\right) } + A  {|\lambda_*|^{N-1} \over \mu_*^{2(N-1)} }  C e^{D {N \over 2^N}\left(i-2^{N-1}\right) } \\
    &\le {|\lambda_{*}|^{N-1} \over \mu_*^{2(N-1)}} \left( |\lambda_*|^{N-1} B^2+A \right)  C e^{-D {N+3 \over 8}} e^{ D {N+1 \over 2^{N+1}} i }.
  \end{align*}

  To prove the bound for the orbits of length between $2^{N-1}+2^{N-2}$ and
  $2^{N} - 1$, one can consider $i=2^{N}-k$, $1\leq k \leq 2^{N-2}$, and repeat
  the calculations closely following second part of Lemma
  \ref{expansion_bounds}.

  The conditions $ {|\lambda_{*}|^{N-1} \mu_*^{2(1-N)}} (|\lambda^{N-1}_{*}|B^2+A) e^{-D
    {N+3 \over 8}} \le 1$ are implied by \eqref{Dcond1}. The claim follows.
\end{proof}

\section{Bounded geometry}

{We will demonstrate existence of unbounded geometry for the
  fixed point map $F^{*}$. Since the Rigidity Theorem~\ref{rigidity} tells us
  that the dynamics of any $F \in W^{s}_{loc}(F_*)$ on their Cantor set is
  conjugate to that of $F_{*}$ by a $C^{1+\alpha}$ transformation, identical
  results hold for all maps in the strong universality class of $F_{*}$. We will
  therefore also use the shorthand $\psi_{i}$ for the rescaling functions of
  $F_{*}$ since $\psi^{n}_{i} = \psi^{m}_{i}$ for all $n,m$ by the fact that
  $R^{n}F_{*} = R^{m}F_{*}$. Note also that because of this the renormalization
  microscope $\Psi_{w} \equiv \Psi_{w}^{F_*}$ is given just by compositions of
  $\psi_{0}$ and $\psi_{1}$. In particular $\Psi_{0^{n}}$ consists of $n$
  compositions of $\psi_{0}$.}
{
  \begin{definition}
    We will say that an infinitely renormalizable Hénon-like map has $D$-bounded
    geometry if $\diam(B_{wi}^n) \asymp_D \operatorname{dist}(B_{w1}^n, B_{w0}^n
    )$, $i \in \{0,1\}$, $w \in \{0,1\}^n$, $n \in \N$, where $\asymp_D$ stands
    for commensurability with a constant $D$. If no such $D$ exists, then we say
    that the geometry is unbounded.
  \end{definition}
}

Before we proceed to the proof of existence of bounded geometry, we will also
mention the following version of the Mean Value Theorem for vector-valued
functions in $\R^{n}$ which will be necessary for both bounded and unbounded
geometry.

\begin{lemma}\label{MVT}
  Let $\UU \subset \R^n$ be open, $F : \UU \mapsto \R^m$ continuously
  differentiable, and $x \in \UU$, $h \in \R^n$ be vectors such that the line
  segment $x + t h$, $0\le t\le 1$ remains in $\UU$. Then we have:
  \[
    F(x+h)-F(x)=\left(\int_{0}^{1} DF (x+th) \ dt\right) \cdot h.
  \]
\end{lemma}

The bounds obtained in Lemmas~\ref{expansion_bounds}, \ref{contraction_bounds}
and \ref{expansion_bounds_second}, together with the Mean Value Theorem, will be
used to demonstrate that pieces in the orbit of certain length of the two
central pieces $B^{n+m}_{0^{n+m}0}$ and $B^{n+m}_{0^{n+m}1}$ have bounded
geometry. The total measure of these iterates is nonzero.

Throughout the proof $C$ will be a stand in for a positive constant whose value
is immaterial to the proof.

\begin{thmA}
  There exists $m \in \N$ such that for all $F \in
  W^{s}_{\mathrm{loc}}(F_{*})$, the measure of pieces $B^{n+m }_{w} \in \BB^{n+m
  }$, $n \in \N$, with bounded geometry is at least $2^{-m-1}$: specifically,
  there exists a constant $C>1$, such that for $v\in\{0,1\}$ and all $w \in
  \{0,1\}^{n}$,
  \begin{equation}
    \label{bounded} C^{-1} < \frac{\dist(B^{n+m }_{w 0^{m} 0}, B^{n+m}_{w 0^{m}1})}{\diam(B^{n+m}_{w 0^{m} v})} < C.
  \end{equation}
\end{thmA}
\begin{remark} \label{quantifiers}
  The pieces appearing in the inequality $(\ref{bounded})$ are iterates up to length $2^n$ of a pair of pieces contained in a deep central piece  $B^{n+m }_{0^{n+m+1}}$. Thus, the quantifier $n+m$ specifies the depth of this pair, while the integer $n$ specifies the length of the orbit along which the bounded geometry has not been destroyed yet by dynamics. 

In the statement of the theorem, both the depth and the length of the orbit can be arbitrarily large by taking $n \rightarrow \infty$, while the ``separation'' $m$ between them is kept constant.

  Finally, the integer $k$, appearing in the proofs below, serves for counting iterates in the orbit of length $2^n$.
\end{remark}
\begin{proof}
  We will prove the bounds for $F_{*}$: the claim itself follows by rigidity in
  $W^{s}_{\mathrm{loc}}(F_{*})$.
  
  Since pieces $B^{n+m }_{0^{n+m }0}$ and $B^{n+m }_{0^{n+m }1}$ are linear
  rescalings of the horizontally separated pieces $B^{0}_{0}$ and $B^{0}_{1}$,
  there exists a constant $K$, such that
  \begin{equation}
    \label{B0bounded}  K^{-1} \le  \frac{\dist(B^{n+m}_{0^{n+m} 0}, B^{n+m}_{0^{n+m} 1})}{\diam(B^{n+m}_{0^{n+m} v})}  \le K
  \end{equation}
  for all natural numbers $n$ and $m$. We consider the action of the iterates
  $F_*^{k}$ on the sets $B^{n+m }_{0^{n+m }v} \subset B^{n+m-1}_{0^{n+m }}$,
  $v=0,1$.
  
  \vspace{3mm}
  \noindent \underline{\it Part 1) An upper bound on the geometry.}
  \vspace{3mm}
  
  We recall that the the central piece $B_0$ can be chosen to be an ellipse, see Remark $\ref{B0ellipse}$, and so is its rescaling $B^{n+m-1}_{0^{n+m}}$.  Consider a horizontal section of the set $B^{n+m-1}_{0^{n+m}}$, and let $c$
  and $d$ be some intersection points of the horizontal line with $\partial
  B^{n+m }_{0^{n+m }0}$ and $\partial B^{n+m }_{0^{n+m }1}$ respectively, such
  that the segment $(c,d)$ does not intersect $B^{n+m }_{0^{n+m}v}$, $v=0,1$. We
  have for the two pieces $B^{n+m }_{w0^{m}v}=F_*^{k} ( B^{n+m }_{0^{m+n}v})$, $k=k(w)$,
  \begin{equation}
     \label{ratio2} \dist \hspace{-1mm} \left(B^{n+m }_{w0^{m }0},B^{n+m }_{w0^{m }1}  \right)   \le   \|D F^k_*(y(k)) e \|  \|c -d\|,
  \end{equation}
  where $e$ is the horizontal unit vector, and the point $y(k) \in (c,d)$ depends on $k=k(w)$. This bound is an immediate consequence of Lemma $\ref{MVT}$.

  By Lemmas~\ref{contraction_bounds} and \ref{expansion_bounds_second} there
  exist $C>0$ and $D>0$, such that
  \begin{align*}
    \| D F^k_* \arrowvert_{B_{0^{n+m}v}^{n+m}} v \|  & \ge C^{-1} e^{-D {n+m \over 2^{n+m}} k},\\
    \| D^2 F^k_* \arrowvert_{B_{0^{n+m}v}^{n+m}}\|  & \le C e^{D {n+m \over 2^{n+m}} k},
  \end{align*}
  for any unit $v$ and for all $1 \leq k < 2^n$, whenever $m \ge 1$. Therefore,
  there exists $\alpha_m>1$, which approaches $1$ as $m \rightarrow \infty$,
  such that
  \begin{align}
    \label{der1}  \| D F^{k}_{*} \arrowvert_{B_{0^{n+m}}^{n+m-1}} \| & \ge C^{-1}\alpha_m^{-n},\\
    \label{sec_der1}  \| D^{2} F^{k}_{*} \arrowvert_{B_{0^{n+m}}^{n+m-1}} \| & \le C\alpha_m^{n}.
  \end{align}
  Bounds \ref{der1} and \ref{sec_der1} imply that the distortion of the
  derivative in a piece is small: For any $x$ and $y$ in $B_{0^{n+m}}^{n+m-1}$,
  and a unit vector $v$
  \begin{equation}
    \label{distort}{\| D F_*^k(y) v\|  \over \| D F_*^k(x) v\|} \le {\| D F_*^k(x) v\| +  C \| D^2 F^k_* \arrowvert_{B_{0^{n+m}}^{n+m-1}}  \| \theta_2^{n+m}  \over \| D F_*^k(x) v\|}   \le 1+C \alpha_m^{2 n} \theta_2^{n+m},
  \end{equation}
  where, if $m$ is sufficiently large, $\alpha_m^{2 n} \theta_2^{n+m} < \beta^n$
  for some $\beta <1$ and any $n \in \N$. Here, $\theta_2$ is an upper bound on the contraction rate of the renormalization microscope, see Proposition $\ref{esti}$.

  Similarly, for a bound from below.
  Additionally,
  \begin{equation}
    \label{distort2}{\| D F_*^k(y) v-  D F_*^k(x) v \|  \over \| D F_*^k(x) v\|}  \le \beta^{n},
  \end{equation}

  We can, therefore, evaluate derivatives in \ref{ratio2} at the tip $\tau$:
  \begin{align}
    \label{ratio22} \dist \left(B^{n+m }_{w0^{m }0},B^{n+m }_{w0^{m }1}  \right)   &\le   (1+C \beta^n) \|D F^k_*(\tau) e \|  \|c -d\|.
  \end{align}

  Similarly, fix $v \in \{0,1\}$ and consider a section of the sets $B^{n+m
  }_{0^{n+m }v}$ by a horizontal line. Let $a$ and $b$ be some intersection
  points in this section with the boundary $\partial B^{n+m}_{0^{n+m }v}$, such
  that the horizontal segments $(a,b) \subset B^{n+m }_{0^{n+m }v}$. Let
  $k(w)=\sum_{i=0}^{n-1} w_{i } 2^{i}$. We have then for $B^{n+m}_{w0^{m}v}=F_*^{k} ( B^{n+m }_{0^{m+n}v})$, using the fact that the minimum of the action of the derivative  $D F_*^k(y)$ at some point $y \in (a,b)$ on the unit horizontal vector $e$ differs from this action at the tip  $\tau \in  B^{n+m-1}_{0^{n+m }}$ by the a bound on the second derivative times the size of $B^{n+m-1}_{0^{n+m }}$: 
  \eqref{ratio22}:
  \begin{align}
\nonumber    {\rm diam} (B^{n+m }_{w0^{m }v})  & \ge \| F_*^k(a)-F_*^k(b) \| \ge \\
\nonumber    & \ge \min_{y \in (a,b)} \| D F_*^k(y) e \| \|a-b\| \ge \\
\nonumber    & \ge \left(\| D F_*^k(\tau) e \|-C \| D^2 F^k_* \arrowvert_{B_{0^{n+m}}^{n+m-1}}  \| \theta_2^{n+m} \right) \|a-b\| \ge \\
 \label{ratio1}    & \ge (1-C \beta^n) \left\| D F^k_*(\tau) e \right\| \|a-b\|,
  \end{align}
  where $\beta$ is the same constant as in $(\ref{distort2})$,  $e=[1,0]$. 
 
  Since $\|a-b\|$ and $\|c-d\|$ are commensurate, so are ${\rm diam}(B^{n+m
  }_{w0^{m }v})$ and $ \dist \left(B^{n+m }_{w0^{m }0},B^{n+m }_{w0^{m }1}
  \right)$.

  \vspace{3mm}
  \noindent \underline{\it Part 2) A lower bound on the geometry.}
  \vspace{3mm}

  To obtain  bounds on the diameters of the pieces $B^{n+m }_{w 0^{m }v}$ and distances between them, we first notice that the initial pieces $B^{n+m }_{0^{n+m}v}$ are very flat: these are deep pieces obtained from $B_v$ by  a diagonal rescaling $\Lambda_*^{n+m}$, which is very incommensurate in the vertical and horizontal dimensions, specifically, the horizontal size of these pieces is of order $|\lambda_*|^{n+m}$ while the vertical is or order $\mu_*^{n+m}$.  We can think of these two pieces as contained in very thin horizontal rectangles, while line segments connecting any two points $x \in B^{n+m }_{0^{n+m}0}$ and $y \in B^{n+m}_{0^{n+m}1}$  in these two pieces lies in very thin horizontal cones.

  To bound $\dist \left( F^{k}_* ( B^{n+m }_{0^{n+m }0}), F^{k}_* ( B^{n+m
    }_{0^{n+m }1}) \right)$ from below we first consider the set of all
  intervals connecting points in the two sets $B^{n+m }_{0^{n+m }0}$ and $B^{n+m
  }_{0^{n+m}1}$. The vectors along which these intervals are oriented lie in
  some closed cone $\bar \Theta_{[1,0]}(A \gamma^{m+n })$,
  $\gamma=\mu_*/|\lambda_*|$, for some constant $A$. Similarly to the derivation of $(\ref{ratio1})$, we can now use bounds $(\ref{der1})$ and $(\ref{sec_der1})$ to estimate the distortion of the first derivative at two points $y$ and $\tau$ in the piece  $ B^{n+m-1}_{0^{n+m }}$, to obtain  
  \begin{align*}
    \dist  & \left( F^{k}_*  ( B^{n+m }_{0^{n+m }0}),  F^{k}_* ( B^{n+m }_{0^{n+m }1})  \right)   \ge  \\
    & \ge \min_{{ v \in \bar \Theta_{[1,0]}(A \gamma^{m+n }), \atop y \in  B^{n+m-1}_{0^{n+m}}}} \hspace{-2mm}  \|D F^k_*(y) \arrowvert_{{\rm span}\{v\} } \|  \dist  \left(B^{n+m }_{0^{n+m }0}, B^{n+m }_{0^{n+m }1}\right) \ge \\
     & \ge C (1-C \beta^n) \min_{{ v \in \bar \Theta_{[1,0]}(A \gamma^{m+n })}} \hspace{-2mm}  \|D F^k_*(\tau) \arrowvert_{{\rm span}\{v\} } \|  \dist  \left(B^{n+m }_{0^{n+m }0}, B^{n+m }_{0^{n+m }1}\right).
  \end{align*}

  To obtain a bound on the diameters of the pieces $B^{n+m }_{0^{n+m }v}$ from above, we first notice that each of the pieces is contained in some rectangle $\RR_v$  with length $2 a_v |\lambda_*|^{n+m}$ and height $2 b_v \mu_*^{n+m}$. We emphasize that $2 a_v |\lambda_*|^{n+m} \gg 2 b_v \mu_*^{n+m}$.

  We claim that the largest possible bound on the diameter of pieces $F^{k}_* ( B^{n+m }_{0^{n+m }v})$ is obtained by the largest possible expansion by $F^k_*$ in a certain thin horizontal cone. Indeed, let $\zeta=(a_v |\lambda_*|^{n+m}, s b_v \mu_*^{n+m})$, $-1 \le s \le 1$, be an interval connecting the center (identified with $(0,0)$) of the rectangle $\RR_v$ with some point on the right vertical side of $\RR_v$, and let  $\xi=(t a_v |\lambda_*|^{n+m},  b_v \mu_*^{n+m})$, $-1 <-\delta < t < \delta < 1$, be an interval connecting the center with a point in a closed subinterval  of  the upper horizontal side of $\RR_v$. Let $D F^k(0,0)$ be of the form
  \begin{equation}
    \label{der_matrix}    D F^k(0,0)=\begin{bmatrix}
    A_k & B_k\\
    C_k  & D_k
    \end{bmatrix}.
  \end{equation}
  Then, fixing $\| \cdot \|$ as the Euclidean distance,  
  \begin{align}
    \nonumber   \|  D F^k_*(0,0) \zeta \|& = \sqrt{ (A_k a_v |\lambda_*|^{n+m} \pp s B_k b_v\mu_*^{n+m})^2 \pp (C_k a_v |\lambda_*|^{n+m}\pp s D_k b_v\mu_*^{n+m})^2} \\
    \label{eqder1}    & = |\mu_*|^{n+m} \sqrt{B_k^2 + D_k^2}  \sqrt{T_1+T_2+T_3}, \\
    \nonumber   \|  D F^k_*(0,0) \xi \|& = \sqrt{ (A_k t a_v |\lambda_*|^{n+m} \pp  B_k b_v\mu_*^{n+m})^2 \pp (C_k t a_v |\lambda_*|^{n+m}\pp D_k b_v\mu_*^{n+m})^2} \\
    \label{eqder2}    & = |\mu_*|^{n+m} \sqrt{B_k^2 + D_k^2}  \sqrt{S_1+S_2+S_3},
  \end{align}
  where
  \begin{align*}
    T_1&={A_k^2+C_k^2 \over (B_k^2 \pp D_k^2) \gamma^{2(n+m)}}  a_v^2,  \hspace{3mm} T_2=2 s a_v b_v  {A_k B_k+C_k D_k \over (B_k^2 \pp D_k^2) \gamma^{n+m} }, \hspace{2mm} T_3=s^2 b_v^2 , \\
    S_1&=t^2 {A_k^2+C_k^2 \over (B_k^2 \pp D_k^2) \gamma^{2(n+m)}}  a_v^2, \hspace{2mm}   S_2=2 t a_v b_v  {A_k B_k+C_k D_k \over (B_k^2 \pp D_k^2) \gamma^{n+m} }, \hspace{3mm} S_3=b_v^2,
  \end{align*}
  Notice, that
  \begin{equation}
    \label{der_ratio}{A_k^2+C_k^2\over B_k^2+D_k^2}={\left\| D F^k(0,0)\left[{1 \atop 0 }\right] \right\|^2 \over \left\| D F^k(0,0)\left[{0 \atop 1 }\right] \right\|^2}.
  \end{equation}
  By Lemmas~\ref{contraction_bounds} and \ref{expansion_bounds} the RHS of $(\ref{der_ratio})$ is bounded from below by $C \alpha_m^{-4 n}$ and from above by $C \alpha_m^{4 n}$, where $\alpha_m$ is as in $(\ref{der1})$. Therefore, we have that 
  \begin{align*}
    {A_k^2+C_k^2\over (B_k^2+D_k^2)} & \ge  C \alpha_m^{-4 n} \gamma^{-2(n+m)}, \\
    {A_k B_k+C_k D_k\over B_k^2+D_k^2} & \le {\sqrt{(A_k^2+C_k^2)(B_k^2+D_k^2)} \over  B_k^2+D_k^2}  \\
    & \le {\sqrt{A_k^2+C_k^2} \over  \sqrt{B_k^2+D_k^2} }  \\
    & \le C \alpha_m^{2 n}.
  \end{align*}
  Since $\alpha_m \rightarrow 1$ as $m \rightarrow \infty$, while $\gamma$ is strictly less than $1$,  the first term $T_1$ in $(\ref{eqder1})$ is exponentially large in $m$ for any fixed $n \in \N$, and $T_2/T_1 \rightarrow 0$ and $T_3/T_1 \rightarrow 0$ as $m \rightarrow \infty$. Similarly for the terms $S_i$.  Therefore, in both equalities $(\ref{eqder1})$ and $(\ref{eqder2})$ the first term under the square root is dominant, and, since, $|t|$ is strictly less than $1$, there exists an $m \in \Z$, such that for all $n \in \NN$ and $1 \le k<2^n$,
  \begin{equation}
    \label{ineq_ders}    \|  D F^k_*(0,0) \zeta \| >  \|  D F^k_*(0,0) \xi \|.
  \end{equation}

 Since the expansion of distances by $F_*^k$ differs from the linear expansion by $D F^k_*$ by a factor $1+C \beta^n$, we conclude that to obtain an upper bound on the diameter,  it is sufficient to consider the expansion in the horizontal cone  $\bar \Theta_{[1,0]}(R \gamma^{m+n })$, for some $R \ge b_v/\delta a_v$. 
  
Therefore,
  \begin{align*}
    {\rm diam} \  F^{k}_* ( B^{n+m }_{0^{n+m }v}) & \le  C \ (1+C\beta^n) \hspace{-4mm} \max_{ v \in \bar \Theta_{[1,0]}(R \gamma^{m+n }) }   \hspace{-4mm} \|D F^k_*(\tau) \arrowvert_{{\rm span}\{v\} }  \| \times \\
                                                  &\hspace{20mm} \times {\rm diam} \  (B^{n+m }_{0^{n+m }v}).
  \end{align*}
  Suppose that the maximum of the norm of the
  derivative above is achieved on some unit  vector $\xi^k$. We get
  \begin{align*}
    \dist   \left( B^{n+m }_{w0^{m }0},  B^{n+m }_{w 0^{m }1} \right)   &\ge   \|D F^k_*(\tau) \zeta^k \| \  c \ \dist\left( B^{n+m }_{0^{n+m }0}, B^{n+m }_{0^{n+m }1} \right), \\
    {\rm diam} (B^{n+m }_{w0^{m }v})  & \le \|D F^k_*(\tau) \xi^k \| \ C \ {\rm diam} \  (B^{n+m }_{0^{n+m }v}),
  \end{align*}
  for some constants $c$ and $C$. Consider the ratio $\|D
  F^k_*(\tau) \zeta^k \| / \|D F^k_*(\tau) \xi^k \|$: Since $\zeta^k$ and
  $\xi^k$ belong to a cone with the opening $\max\{A,R\} \gamma^{n+m}$,
\comment{  
  
  \textcolor{red}{
    \begin{align*}
      {\|D F^k_*(\tau) \zeta^k \| \over  \|D F^k_*(\tau) \xi^k \|}  & \le   {\|D F^k_*(\tau) \xi^k \| + \|D F^k_*(\tau) \| \|\zeta^k-\xi^k\|  \over  \|D F^k_*(\tau) \xi^k \|} \\
                                                                    & \le   {\|D F^k_*(\tau) \xi^k \| + C \alpha_m^n \gamma^{n+m}  \over  \|D F^k_*(\tau) \xi^k \|} \\
                                                                    & \le  1 + C \alpha_m^{2n} \gamma^{n+m} \\
                                                                    & \le  1 + C \beta^n
  \end{align*}}
  }
{
    \begin{align*}
      {\|D F^k_*(\tau) \zeta^k \| \over  \|D F^k_*(\tau) \xi^k \|}  & \geq {\|D F^k_*(\tau) \zeta^k \| \over \|D F^k_*(\tau) \zeta^k \| + \| DF^{k}_{*}(\tau)\| \|\xi^{k} - \zeta^{k}\|} \\
                                                                    & \geq {\|D F^k_*(\tau) \zeta^k \| \over \|D F^k_*(\tau) \zeta^k \| + C \alpha_m^n \gamma^{n+m}} \\
                                                                    & \geq \frac{1}{1 + C \alpha_m^{2n} \gamma^{n+m}} \\
                                                                    & \geq \frac{1}{1 + C \beta^n}
    \end{align*}}
  for some $C>0$ and $\beta<1$.

  The lower bound follows.
\end{proof}

\section{Unbounded geometry}
\subsection{Unbounded geometry near the tip}

As in the case of  bounded geometry, we will show existence of positive
  measure of unbounded geometry for the fixed point $F_{*}$. The result for any
  $F\in \Wl$ follows by rigidity.
\comment{ 
  \textcolor{red}{We will demonstrate existence
  of unbounded geometry for the fixed point map $F_*$. Since the Rigidity
  Theorem~\ref{rigidity} tells us that the dynamics of any $F \in W^{s}(F_*)$ on
  their Cantor set is conjugate to that of $F_*$ by a $C^{1+\alpha}$
  transformation, identical results hold for all maps in the strong universality
  class of $F_*$. We will therefore also use the shorthand $\psi_{i}$ for the
  rescaling functions of $F_{*}$ since $\psi^{n}_{i} = \psi^{m}_{i}$ for all
  $n,m$ by the fact that $R^{n}F_{*} = R^{m}F_{*}$. Note also that because of
  this the renormalization microscope $\Psi_{w} \equiv \Psi_{w}^{F_*}$ is given
  just by compositions of $\psi_{0}$ and $\psi_{1}$. In particular
  $\Psi_{0^{n}}$ consists of $n$ compositions of $\psi_{0}$.}

\textcolor{red}{
  \begin{definition}
    We will say that an infinitely renormalizable Hénon-like map has $D$-bounded
    geometry if $\diam(B_{wi}^n) \asymp_D d(B_{w1}^n, B_{w0}^n )$, $i \in
    \{0,1\}$, $w \in \{0,1\}^n$, $n \in \N$, where $\asymp_D$ stands for
    commensurability with a constant $D$. If no such $D$ exists, then we say
    that the geometry is unbounded.
  \end{definition}
}
}

We begin with the following simple observation. Recall that the tip
$\tau=(0,0)$ is the fixed point of $\psi_{0}$. One can also find the fixed point
of the second map, $\psi_{1}$, in the iterated function system.

\begin{lemma}
  \label{fp2}
  $T(F_{*}(\tau)) \in B_{1}$ is the fixed point of $\psi_{1}$
\end{lemma}
\begin{proof}
  We have according to~\eqref{sdef}
  \[
    F_{*}(\tau)=\left( 1 \atop s^{*}(0,1) \right) \in B_1.
  \]
  and $T(F_{*}(\tau))=(1,-s^{*}(0,1))$. According to the midpoint
  equation~\eqref{midpoint}, together with the normalization $z(1,0)=z(0,1)=1$,
  \[
    -s^{*}(0,1)=s^{*}(\lambda_*,1),
  \]
  and, therefore, according to~\eqref{sdef}, the preimage of
  $T(F_{*}(\tau))=(1,s^{*}(\lambda_*,1))$ under $F_*$ is
  $(\lambda_*,-s^{*}(1,\lambda_*))$. By the fixed point equation~\eqref{Rs} we
  have $-s^{*}(1,\lambda_*)=-\mu_* s^{*}(0,1)$, and the preimage of
  $T(F_{*}(\tau))$ is the point $(\lambda_*, -\mu_* s^{*}(0,1))$.

  We conclude that $(1,-s^{*}(0,1))=T(F_{*}(\tau))$ is the fixed point of the
  contraction $\psi_1$:
  \begin{eqnarray}
    \nonumber T(F_{*}(\tau))&=&F_{*}((\lambda_*, -\mu_* s^{*}(0,1))) \\
    \nonumber &=&F_{*} \circ \psi_{0}(1,-s^{*}(0,1)\\
    \nonumber &=&\psi_1(1,-s^{*}(0,1)) \in B_1.
  \end{eqnarray}
\end{proof}

\begin{figure}
  \centering
      {\includegraphics[scale=0.23]{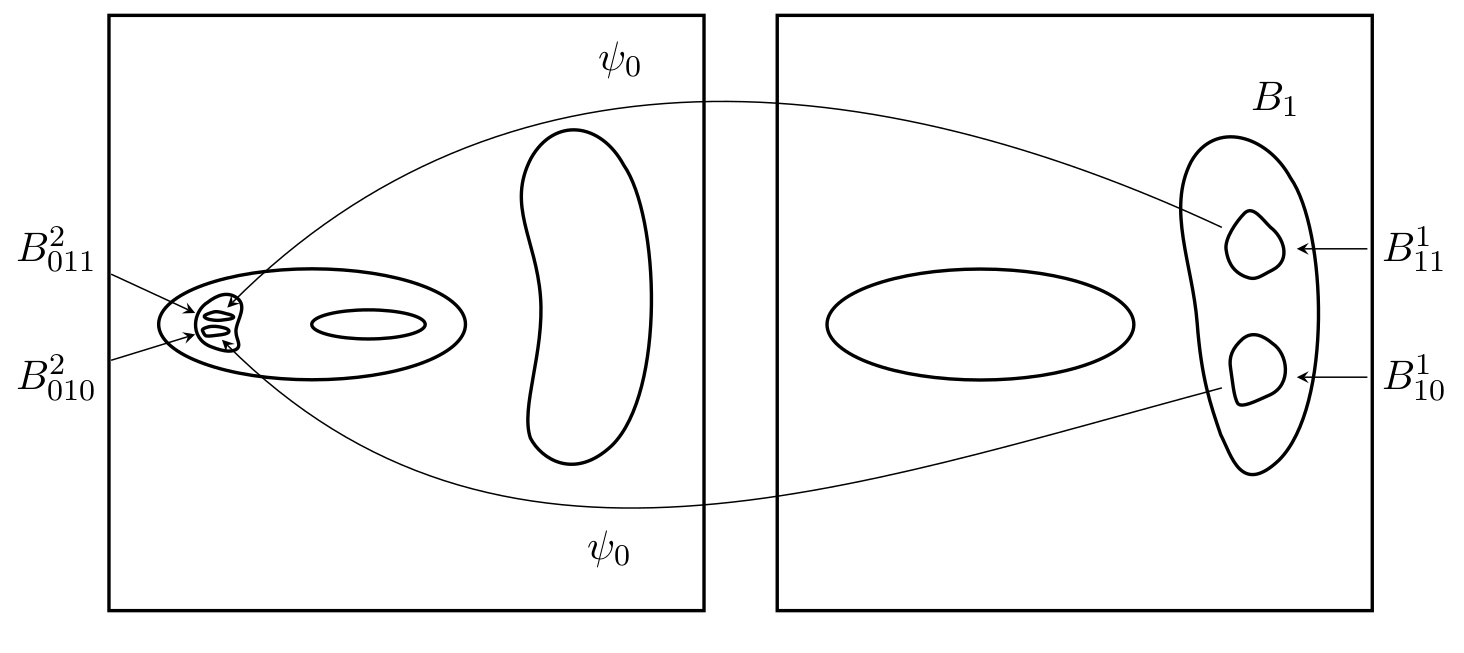}}
\caption{Unbounded geometry near the tip.}
  \label{fig:unbounded_tip}
\end{figure}

We will now describe the geometry of the Cantor set near the tip.

Relative to their sizes the pieces $B_{011}^2$ and $B_{010}^2$ are closer together, thinner and wider than their preimages $B_{11}^1$ and $B_{10}^1$ due to the different scalings of $\psi_0$ in the vertical and horizontal directions. This geometry will only become more extreme at $B_{0^n11}^{n+1}$ and $B_{0^n10}^{n+1}$ as more and more iterates of $\psi_0$ are applied. Figure~\ref{fig:unbounded_tip} illustrates this idea behind the proof of unbounded geometry near the tip.

\begin{proposition}\label{unb_tip}
  The geometry of $\CC_{F_{*}}$ is unbounded near the tip.
\end{proposition}
\begin{proof}
  Consider the two pieces $B_{10}^1$ and $B_{11}^1$. We have,
  using $\psi_{0}^{2}(\tau) = \tau$,
  \[
    B_{10}^1 \ni \psi_1(\psi_0(\tau))=F_{*}(\Lambda_*^2(\tau))=F_*(\tau)=\left(1
      \atop s^{*}(0,1) \right).
  \]
  On the other hand, the fixed point of $\psi_1$,
  $T(F_{*}(\tau))=\psi_1(T(F_{*}(\tau))) \in B_{11}^1$.
  Now, consider the pieces $B_{0^n 1 1}^{n+1}$ and $B_{0^n 1 0}^{n+1}$. We have
  \[
    (\lambda_*^n, -\mu_*^n s^{*}(0,1))=(\lambda_*^n, \mu_*^n s^{*}(\lambda_*,1))
    \in B_{0^n 1 1}^{n+1} \ {\rm and} \ (\lambda_*^n, \mu_*^n s^{*}(0,1)) \in
    B_{0^n 1 0}^{n+1}.
  \]
  Therefore,
  \[
    \dist ( B_{0^n 1 1}^{n+1}, B_{0^n 1 0}^{n+1} ) < 2 \mu_*^n |s^{*}(0,1)| \sim \mu_*^n,
  \]
  while, since $\diam (B_{0^{n}w}^n) \geq \diam (\pi_{x}(B_{0^{n}w}^n)) = \diam
  (\pi_{x}B_{w})\lvert \lambda_{*} \rvert^{n}$,
  \[
    \diam (B_{0^n 1 1}^{n+1}) \asymp_D \diam (B_{0^n 1 0}^{n+1}) ) \geq
    \min\left\{ \diam(\pi_x B_{1 1}^1), \diam(\pi_x B_{1 0}^1) \right\}
    |\lambda_*|^{n},
  \]
  where
  \[
    D=\max\left\{ {\diam(\pi_x B_{1 1}^1) \over \diam ( \pi_x B_{1 0}^1)},
      {\diam(\pi_x B_{1 0}) \over \diam (\pi_x B_{1 1}^1)} \right\}.
  \]
  The conclusion follows.
\end{proof}

\subsection{Positive measure of unbounded geometry}



To prove positive measure of unbounded geometry we will use the unbounded
geometry near the tip together with the techniques used to prove bounded
geometry in order to spread the unbounded geometry through part of $\BB^{n+m}$
using $F_{*}$. Rigidity then shows that the same is true for any $F\in
W^{s}_{\mathrm{loc}}(F_{*})$.

We begin by noticing that $B_{0^{n+m} 1 v}^{{n+m}+1} \subset B_{0^{n+m}}^{{n+m}-1}$ for
 $v \in \{0,1\}$.   Recall that the the central piece $B_0$ has been constructed as an ellipse, see Remark $\ref{B0ellipse}$, therefore, $B^{{n+m}-1}_{0^{{n+m}}}$ is an ellipse as well, and we can choose vertical line segments $(c^{{n+m}}, d^{{n+m}})\subset B^{{n+m}-1}_{0^{{n+m}}}$ with endpoints $c^{{n+m}}\coloneqq \psi_0^{{n+m}}(F_{*}(\tau))$ and $d^{{n+m}}\coloneqq \psi_0^{{n+m}}(T(F_{*}(\tau)))$. Since by Proposition $\ref{unb_tip}$, $c^{{n+m}} \in B_{0^{n+m} 1 0}^{{n+m}+1}$  and $d^{{n+m}}\in B^{{n+m}-1}_{0^{{n+m}}11}$, we will use upper estimates on the expansion of the segment $(c^{{n+m}}, d^{{n+m}})$ by dynamics as an upper bound on the distance between the sets $B_{0^{n+m} 1 0}^{{n+m}+1}$ and $B_{0^{n+m} 1 1}^{{n+m}+1}$.

Similarly to the proof of Theorem A, we will use below two integers $m$ and $n$, and a counting index $k$. The quantifiers $m$ and $n$ will play the same role as before, as explained in Remark $\ref{quantifiers}$: $m+n$ will specify the depth of the pieces, while $n$ specifies the length of the orbit considered. Eventually, $m$ will be a fixed large integer, while $n$ is any natural number. Index $k=k(w)$ counts the iterates of the map $F_*$ in the orbit of length $2^n$.

Let $m \geq 0$. By the Mean Value Theorem~\ref{MVT} we get
\begin{align*}
  \dist( F^{k}_{*}(B_{0^{n+m} 1 0}^{n+m+1}),  F^{k}_{*}(B_{0^{n+m} 1 1}^{n+m+1})  ) & \leq \|{F_{*}^{k}(c^{n+m}) - F_{*}^{k}(d^{n+m})}\| \\
                                                                & \leq \|DF_{*}^{k}\vert_{B^{n+m-1}_{0^{n+m}}}\|\|c^{n+m} - d^{n+m}\| \\
                                                                & \leq 2\mu_{*}^{n+m}\|DF_{*}^{k}\vert_{B^{n+m-1}_{0^{n+m}}}\|\vert s^{*}(0,1)\vert.
\end{align*}

Next we apply Lemma~\ref{expansion_bounds} to bound the influence of
$\|DF_*^{k}\vert_{B^{n+m-1}_{0^{n+m}}}\|$ for $k < 2^{n+m-1}$ and find
\begin{equation}
  \label{dist_i}
  \dist( F_{*}^{k}(B_{0^{n+m} 1 1}^{n+m+1}),  F_{*}^{k}(B_{0^{n+m} 1 0}^{n+m+1}) ) < 2 \mu_{*}^{n+m} \vert s^{*}(0,1)\vert C e^{D {n+m \over 2^{n+m}} k};
\end{equation}
for $0\leq k < 2^{n+m-1}$.

Similarly, pick horizontal line segments $(c^{{n+m}}_{v}, d^{{n+m}}_{v})\subset B^{{n+m}+1}_{0^{{n+m}} 1 v}$ by first picking $c^{0}_{v}\in \partial B^{1}_{1v}$ and
$d^{0}_{v}\in \partial B^{1}_{1v}$ such that the horizontal line segment
$(c^{0}_{v}, d^{0}_{v})\subset B^{1}_{1v}$ and letting $(c^{{n+m}}_{v}, d^{{n+m}}_{v}) =\psi_0^{{n+m}}((c^{0}_{v}, d^{0}_{v}))$. Note that $\| c^{{n+m}}_{v} - d^{{n+m}}_{v}\| = \vert\lambda_{*}\vert^{{n+m}}\|c^{0}_{v} - d^{0}_{v}\|$. For
notational convenience we will use
\begin{equation*}
  K^{-1}\coloneqq \min_{v\in\{0,1\}}\|c^{0}_{v} - d^{0}_{v}\|.
\end{equation*}
Following the argument leading to \eqref{ratio22} we get, for $m$ large enough, $n \in \NN$ and $k < 2^{n-1}$,
\begin{align*}
  \diam(F^{k}_{*}(B^{n+m+1}_{0^{n+m}1v})) & \geq \| F^{k}_{*}(c^{n+m}_{v}) - F^{k}_{*}(d^{n+m}_{w}) \| \\
  & \geq (1 - C\beta^n) \min_{ u \in \R^{2}\setminus \{0\}}\frac{\|DF_*^{k}(\tau)u \|}{\|u\|} \|  \|c^{n+n}_{v} - d^{n+m}_{v}\| \\
  & \geq K^{-1}\vert\lambda_{*}\vert^{n+m}(1 - C\beta^{n})\min_{u\in \R^{2} \setminus \{0\} }\frac{\|DF_*^{k}(\tau)u\|}{\|u\|}
\end{align*}
where $C$ is some constant and $\beta < 1$. We have used again that the expansion by $F^k_*$, $k<2^{n-1}$, in deeps pieces coincides with the expansion by the derivative $D F^k_*(\tau)$ at the tip up to a factor $1-C \beta^n$. 

Now we can use Lemma~\ref{contraction_bounds} to bound the influence of
$\min_{u\in \R^{2} \setminus \{0\}}\frac{\|DF_*^{k}(\tau)u\|}{\|u\|}$:
\begin{equation}
  \label{diam_i}  \diam (F^{k}_{*}(B_{0^{n+m} 1 v}^{n+m+1})) > K^{-1} |\lambda_{*}|^{n+m}(1 - C \beta^{n}) C^{-1} e^{-D {n+m \over 2^{n+m}} k}, 
\end{equation}
where $v \in \{0,1\}$ and $D$ verifies the hypothesis of Lemmas~\ref{expansion_bounds} and~\ref{contraction_bounds} and the constant $C$, for simplicity, can be chosen to be the same in both \eqref{dist_i} and \eqref{diam_i}.

With these two bounds we now see that 
\begin{align*}
  \frac{\dist(F_{*}^{k}(B^{n+m+1}_{0^{n+m}10}), F_{*}^{k}(B^{n+m+1}_{0^{n+m}11}))}{\diam(F_{*}^{k}(B^{n+m+1}_{0^{n+m}1v}))} & \leq\frac{2 \mu_{*}^{n+m} \vert s^{*}(0,1)\vert C e^{D {n+m \over 2^{n+m}} k}}{K^{-1} |\lambda_{*}|^{n+m}(1 - C \beta^{n})  C^{-1} e^{-D {n+m \over 2^{n+m}} k}} \\
  & = 2KC^{2}\vert s^{*}(0,1)\vert \frac{\mu^{n+m}_{*}}{\vert\lambda_{*}\vert^{n+m}}\frac{1}{1 - C \beta^{n}}e^{D\frac{n+m}{2^{n+m-1}}k}
\end{align*}
for all $k < 2^{n-1}$ and  $v \in \{0,1\}$. We see that this will approach $0$ as $n\to\infty$ if
\begin{equation*}
  (n+m)(\ln\mu_{*} - \ln\vert\lambda_{*}\vert) + D\frac{n+m}{2^{n+m-1}}k < 0
\end{equation*}
or, equivalently,
\begin{equation*}
  D < \frac{2^{n+m-1}}{k}(\ln\vert\lambda_{*}\vert - \ln\mu_{*}).
\end{equation*}
Since $k < 2^{n-1}$ we see that this will be true for all $0 \leq k < 2^{n-1}$
if
\begin{equation}
  \label{eq:D_condition_bounded_geometry}
  D < 2^{m}(\ln\vert\lambda_{*}\vert - \ln\mu_{*}),
\end{equation}
which can be ensured by choosing $m$ sufficiently large. We notice that the measure of the pieces
\[
  \bigcup_{k=0}^{2^{n-1}-1} \left(F^{k}_{*}(B_{0^{n+m} 1 0}^{n+m+1}) \cup F^{k}_{*}(B_{0^{n+m} 1 1}^{n+m+1}) \right)
\]
in $\BB^{n+m+1}$, is $2^{n-1}2^{-n-m-2}=2^{-m-3}$. Using rigidity we have
arrived at the following result, where introducing an auxiliary parameter
$\gamma$ allows us to control the rate at which the geometry becomes unbounded.


\begin{thmB}
  Let $F \in W^{s}_{\mathrm{loc}}(F_*)$. For every $1\leq \gamma <
  |\lambda_*|\mu^{-1}_{*}$ there exists an $m = m(F,\gamma)\in\N$ such that the
  measure of the pieces in $\BB^{n+m}$ satisfying
  \[
    \diam (B_{w v}^{n+m}) > \gamma^{n+m-1} \dist( B_{w 1}^{n+m}, B_{w 0}^{n+m} ), \
    w \in \{0,1\}^{n+m}, \ v \in \{0,1\},
  \]
  for all $n\in\N$ is at least $2^{-m-4}$.
\end{thmB}

The following result is a straight forward corollary of the above theorem.

\begin{corollary}
  The measure of pieces in $\BB^{n}$ with unbounded geometry is positive.
\end{corollary}

\section{Ergodic properties of the derivative cocycles}\label{ontheaverage}

In this section we will study some ergodic properties of the derivative cocycles
of $F_*$.

The main results of this Section are Propositions
~\ref{zero_average_cone_expansion}, ~\ref{adequate_conv} and~\ref{angles3}. The
first of them shows that the derivative cocycles do not expand or contract
angles between vectors on the average, the second - that the Birkhoff averages
of continuous functions converge quickly, while the last demonstrates that the
distortion of arbitrarily long microscope branches is uniformly bounded.

\subsection{Distortion of vectors by the microscope branches.} We will start
with a simple result about the distortion of vectors by the derivatives of the
renormalization microscope.

We continue to denote the microscope branches of the renormalization fixed point
$F_*$ by $\psi_0$ and $\psi_1$ and $\Psi_w$.

\begin{lemma}\label{further_bounds} $\phantom{a}$\\
  \noindent $1)$ There exists a constant $0<K<1$ such that for any $k \in \{0\}
  \cup \N$, any non-zero vector $u$ and any two $x, y \in B_{w}^k$, $ w \in
  \{0,1\}^{k+1}$, and for $v \in \{0,1\}$
  \begin{equation}
    \label{DPsi_ratio1} \left( 1-K\theta_2^{k+1} \right)^{|v|}   \le{ \left\| D \psi_v(x) u \right\| \over \left\| D\psi_v(y) u \right\| } \le \left( 1+K\theta_2^{k+1} \right)^{|v|},
  \end{equation}
  \begin{equation}
    \label{DPsiInv_ratio1} \left(1-K\theta_2^{k+1}\right)^{|v|} \le { \left\| \left[ D \psi_v(x)\right]^{-1} u \right\| \over \left\| \left[ D \psi_v(y) \right]^{-1} u \right\| } \le \left( 1+K\theta_2^{k+1}\right)^{|v|},
  \end{equation}
  where $\theta_2$ is as in Proposition $\ref{esti}$.

  \noindent $2)$ For any three unit vectors $e$, $h$ and $k$, the following
  holds.
  \begin{align}
    \label{DF_ratio2} {  \left\| D \psi_v(x) e \right|  \over  \left\| D \psi_v(x) h \right\| } &= {\sin(e,k) \over \sin(h,k)     } { \sin \left( D \psi_v(x) h,  D \psi_v(x) k \right)  \over  \sin \left( D \psi_v(x) e,  D \psi_v(x) k \right)      }, \\
    \label{DFInv_ratio2}{  \left\| \left[ D \psi_v(x) \right]^{-1} e \right\|  \over  \left\| \left[ D \psi_v(x)  \right]^{-1} h \right\| }& ={\sin(e,k) \over \sin(h,k)     } { \sin \left( \left[ D \psi_v(x) \right]^{-1} h, \left[ D \psi_v(x) \right]^{-1} k \right)  \over  \sin \left( \left[ D \psi_v(x) \right]^{-1} e, \left[ D \psi_v(x) \right]^{-1} k \right)      }.
  \end{align}
  where $(u,v)$ denotes the angle between vectors $u$ and $v$.
\end{lemma}
\begin{proof} $\phantom{aaa}$

  \medskip

  \noindent $1)$ Claims  $\ref{DPsi_ratio1}$ and $\ref{DPsiInv_ratio1}$ follow by first order approximation of derivatives of the smooth maps $\psi_v$. For any unit vector $u$, and $x$, $y$ contained in some compact subset of the domain of definition of $\psi_v$, such as $B_w^n$, we have $\| D \psi_v(x) u- D \psi_v(y)  u\| \le C \|x-y\|$.    Since $x, y \in B_{w}^k$, and the sizes of these pieces reduce as $\theta_2^k$, we get that  $\|x-y\| \le C \theta_2^k$. By Proposition  $\ref{esti}$,  $\|D \psi_v(y) u\| \ge C \theta_1$, therefore
  \begin{equation}
{\| D \psi_v(x)  u -  D \psi_v(y)  u\| \over  \| D \psi_v(y)  u\|}  \le C \theta_2^k, 
  \end{equation}
which implies the upper bound (recall that $C$ denotes any constant). The lower bound is proved in a similar way.

    Notice, that if $v=0$, then the denominator and   numerator in \eqref{DPsiInv_ratio1} are identical: this explains the power   $|v|$ in the bound.
  
  \medskip

  \noindent $2)$ Notice that the Jacobian of {$\psi_{v}^{-1}(x)$} is constant and equal to
  $|\lambda_*|^{-1} \mu_*^{-1}$. Therefore, for example, to obtain
  \eqref{DFInv_ratio2}, we have
    \begin{align}
     \nonumber  { \left\| \left[D \psi_v(x) \right]^{-1} e \right\|  \over  \left\| \left[  D \psi_v(x) \right]^{-1} h \right\|   }&=   { \left\| \left[  D \psi_v(x) \right]^{-1} e \right\|  \left\|  \left[ D \psi_v(x) \right]^{-1} k \right\| \over   \left\|\left[   D \psi_v(x) \right]^{-1} h \right\|  \left\| \left[  D \psi_v(x) \right]^{-1} k \right\| } \\
     \nonumber &={ \|e\| \|k\| \over \| h\| \|k\| } {\sin(e,k) \over \sin(h,k)     } { \sin \left( \left[ D \psi_v(x) \right]^{-1} h, \left[ D \psi_v(x) \right]^{-1} k \right)  \over  \sin \left( \left[ D \psi_v(x) \right]^{-1} e, \left[ D \psi_v(x) \right]^{-1} k \right)      } \\
    \nonumber &= {\sin(e,k) \over \sin(h,k)     } { \sin \left( \left[ D \psi_v(x) \right]^{-1} h, \left[ D \psi_v(x) \right]^{-1} k \right)  \over  \sin \left( \left[ D \psi_v(x) \right]^{-1} e, \left[ D \psi_v(x) \right]^{-1} k \right)      }.
    \end{align}
    Similarly, for \eqref{DF_ratio2}.
\end{proof}

We will use the following notation in the remainder of the paper: 
\begin{align}
  \label{barrho} \bar \rho_i& =(\rho_1, \rho_2, \ldots, \rho_i),\\
  \label{tilderho} \tilde \rho_{i+1} & =(\rho_{i+1}, \rho_{i+2}, \ldots, \rho_l)
\end{align}
for $i \ge 1$ and any $\rho \in \{0,1\}^l$. Additionally, for any $k \in \N$ and any $\pi \in \{0,1\}^k$, we will denote images of arbitrary points by renormalization microscopes as
\begin{equation}
  \label{xpi} x_\pi = D \Psi_{\pi}  (x),
\end{equation}
while images of indexed points $x_w \in B_w^{n-1}$, with $w \in \{0,1\}^n$, will be denoted by
\begin{equation}
  \label{xpiw} x_{\pi w} = D \Psi_{\pi}  (x_w).
\end{equation}
Similarly, given a vector $u \in \R^2$, denote vectors in its forward and backward orbits as 
\begin{align}
  \label{upi} u_\pi(x) & = D \Psi_{\pi} (x) u, \\
  \label{upim} u_\pi^-(x) & = \left[ D \Psi_{\pi} \right]^{-1}(x) u.
\end{align}
Finally, if the dependence of the vectors in the orbit on the point is not important to the argument, this dependence will be omitted.

We can now evaluate the distortion of vectors by inverse long branches of the renormalization microscope.

\begin{lemma}\label{more_bounds}
 The following holds for any $n,l \in \N$ satisfying $n > l$, any   non-zero vector $u$, any $x_{w}, y_{w} \in B_{w}^{n-1}$, $w\in\{0,1\}^{n}$:
   \begin{equation*}
      \Pi_l \prod_{i=0}^{l-1}  \left( 1  -   K\theta_2^{n\shortminus i} \right)^{  \hspace{-0.5mm} |w_{i+1}|}  \le { \left\| \left[ D \Psi_{\bar w_l}  \right]^{-1} \hspace{-1mm}  (x_w) u \right\| \over \left\| \left[ D\Psi_{\bar w_l}  \right]^{-1} \hspace{-1mm} (y_w)   u \right\| } \le  \Pi_l \prod_{i=0}^{l-1}    \left( 1  +   K\theta_2^{n\shortminus i} \right)^{  \hspace{-0.5mm} |w_{i \shortplus 1}|},
   \end{equation*}
   where  
\begin{equation}
  \nonumber \Pi_1=1, \quad  \Pi_l= {\sin\left(e_{u_{\bar w_1}^{-}(x_w)} ,e_{v_{\bar w_1}^{-}(x_w)} \right) \over \sin\left(e_{u_{\bar w_1}^{-}(y_w)} , e_{v_{\bar w_1}^{-}(x_w)} \right) }  {\sin\left(e_{u_{\bar w_{l}}^{-}(y_w)} ,e_{v_{\bar w_{l}}^{-}(x_w)} \right) \over \sin\left(e_{u_{\bar w_{l}}^{-}(x_w)} , e_{v_{\bar w_{l}}^{-}(x_w)} \right)     }, \quad l \ge 2,
\end{equation}
and
\begin{equation}
  v_{\bar w_i}^{-}(x_w)=\left[ D\Psi_{\bar w_i}  \right]^{-1}(x_w) v,
\end{equation}
is the backward orbit of an arbitrary vector in $\R^2 \setminus \{0\}$, and 
\begin{equation*}
u_{\bar w_i}^{-}(x_w)=\left[ D\Psi_{\bar w_i}  \right]^{-1}(x_w) u, \quad u_{\bar w_i}^{-}(x_w) =\left[ D\Psi_{\bar w_i}  \right]^{-1}(x_w) u
\end{equation*}
with $\bar w_i$ as in $(\ref{barrho})$, $u_{\bar w_0}^{-}=u$ and  $v_{\bar w_0}^{-}=v$.
\end{lemma}
\begin{proof}
For brevity, denote 
\begin{equation}
\label{PT} P_i^{x_w,y_w,u}= { \left\|\left[ D \Psi_{\bar w_i} \right]^{-1}(x_w) u \right\| \over \left\| \left[ D\Psi_{\bar w_i}  \right]^{-1} (y_w) u \right\| }.
\end{equation}
We proceed by induction on $i$. The base of induction for $i=1$ follows from \eqref{DPsiInv_ratio1}, with $k=n-1$ and with $K$ as in \eqref{DPsiInv_ratio1}.  Assume the result for $i \le l$. Consider $x_{\tilde w_{l+1}}, y_{\tilde w_{l+1}} \in B_{\tilde w_{l+1}}^{n+l-1}$.
 \begin{align}
   \nonumber P_{l+1}^{x_w,y_w,u} & =  {  \left\|  \left[ D \Psi_{\bar w_{l+1}} \right]^{-1} (x_w)  u\right\|  \over    \left\|  \left[ D \Psi_{\bar w_{l+1}} \right]^{-1} (y_w) u\right\|    }  \\
   \nonumber & =  {  \left\|  \left[ D \psi_{w_{l+1}}\right]^{-1} (x_{\tilde w_{l+1}})  \left[ D \Psi_{\bar w_{l}} \right]^{-1} (x_w)  u\right\|  \over    \left\|  \left[ D \psi_{w_{l+1}}\right]^{-1} (y_{\tilde w_{l+1}})  \left[ D \Psi_{\bar w_{l}}\right]^{-1} (y_w)  u\right\|} \\
   \nonumber & =  {  \left\|  \left[ \hspace{-0.3mm}   D \psi_{w_{l+1}} \right]^{-1} (x_{\tilde w_{l+1}}) e_{u_{\bar w_l}^{-}(x_w)}\right\|  \over    \left\|  \left[ D \psi_{w_{l+1}} \right]^{-1} (y_{\tilde w_{l+1}}) e_{u_{\bar w_l}^{-}(y_w)} \right\|} P_{l}^{x_w,y_w,u} \\
    \nonumber & =  {  \left\|  \left[   \hspace{-0.3mm}   D \psi_{w_{l+1}}    \hspace{-0.3mm}  \right]^{\shortminus 1}  \hspace{-1mm}  (x_{\tilde w_{l+1}}) e_{u_{\bar w_l}^{-}(x_w)} \right\|  \over    \left\|  \left[    \hspace{-0.3mm}  D \psi_{w_{l+1}}   \hspace{-0.3mm}   \right]^{\shortminus 1}  \hspace{-1mm}   (x_{\tilde w_{l+1}}) e_{u_{\bar w_l}^{-}(y_w)} \right\|} {  \left\|  \left[    \hspace{-0.3mm}  D \psi_{w_{l+1}}   \hspace{-0.3mm}   \right]^{\shortminus 1}  \hspace{-1mm}   (x_{\tilde w_{l+1}}) e_{u_{\bar w_l}^{\shortminus}(y_w)}\right\|  \over    \left\|  \left[  \hspace{-0.3mm}  D \psi_{w_{l+1}}     \hspace{-0.3mm}  \right]^{\shortminus 1}  \hspace{-1mm}   (y_{\tilde w_{l+1}}) e_{u_{\bar w_l}^{-}(y_w)} \right\|} P_{l}^{x_w,y_w,u},
 \end{align}
 where $e_u$ denotes a unit vector in the direction of $u \ne 0$.   We use  $(\ref{DFInv_ratio2})$ on the first ratio  and $(\ref{DPsiInv_ratio1})$ on the second (notice, $x_{\tilde w_{l+1}}$, $y_{\tilde w_{l+1}}$ are in $B_{\tilde w_{l+1} }^{n-l-1}$):
 \begin{align}
   \nonumber P_{l+1}^{x_w,y_w,u} & \le {\sin\left(e_{u_{\bar w_l}^{-}(x_w)} ,e_{v_{\bar w_l}^{-}(x_w)} \right) \over \sin\left(e_{u_{\bar w_l}^{-}(y_w)} , e_{v_{\bar w_l}^{-}(x_w)} \right) }  {\sin\left(e_{u_{\bar w_{l+1}}^{-}(y_w)} ,e_{v_{\bar w_{l+1}}^{-}(x_w)} \right) \over \sin\left(e_{u_{\bar w_{l+1}}^{-}(x_w)} , e_{v_{\bar w_{l+1}}^{-}(x_w)} \right)     }    \times \\
   \nonumber & \hspace{30mm} \times \left(1+ K  \theta_2^{n-l}\right)^{|w_{l+1}|}  P_{l}^{x_w,y_w,u} \\
   \nonumber  & \le \Pi_{l+1}  \prod_{i=0}^{l} \hspace{-0.6mm}  \left( 1  \hspace{-0.5mm} + \hspace{-0.5mm}  K\theta_2^{n - i} \right)^{  \hspace{-0.5mm} |w_{i \shortplus 1}|}.
 \end{align}
The bound from below follows in a similar way.
\end{proof}

The distortion of vectors by direct long branches of the renormalization microscope can be evaluated in an analogous way (we will omit the proof, almost identical to that of Lemma $\ref{more_bounds}$).

\begin{lemma}\label{more_bounds_direct}
 The following holds for any $n,l \in \N$ satisfying $n > l$, any unit vectors $e$ and $k$, any $x, y \in \CC_{F_*}$ and  $w\in\{0,1\}^{n}$:
   \begin{equation*}
      \Pi_l \prod_{i=n}^{l}  \left( 1-  K\theta_2^{n\shortminus i} \right)^{  \hspace{-0.5mm} |w_{i}|}  \le { \left\| D \Psi_{\tilde w_l}   (x) u \right\| \over \left\| D\Psi_{\tilde w_l} (y)   u \right\| } \le  \Pi_l \prod_{i=n}^{l}  \hspace{-0.6mm}   \left( 1  +  K\theta_2^{n\shortminus i} \right)^{  \hspace{-0.5mm} |w_{i}|},
   \end{equation*}
   where  
\begin{equation}
  \nonumber \Pi_n=1, \quad  \Pi_l= {\sin\left(e_{w_n}(x) ,k_{w_n}(x) \right) \over \sin\left(e_{w_n}(y) , k_{w_n}(x) \right) }  {\sin\left(e_{\tilde w_{l}}(y) ,k_{\tilde w_{l}}(x) \right) \over \sin\left(e_{\tilde w_{l}}(x), k_{\tilde w_{l}}(x) \right)     }, \quad l <n,
\end{equation}
where $e_{\tilde w_{l}}(x)$ and $k_{\tilde w_{l}}(x)$ are defined as in $(\ref{upi})$.
\end{lemma}

\subsection{Ergodic properties of the projectivization of the derivative cocycle.}\label{ergodic_prop}

Let $L(x)$ be the projectivization of the action of the derivative $DF$ on
$\R P^1$:
\begin{equation}
  \label{L}  L(x) [v]={ D F(x) v \  \mod \ Z(\R^2)},
\end{equation}
where $[v]$ is the equivalence class of $v$ and $Z(\R^2)$ is a group of non-zero
scalar transformations. Notice that $L(x) \in PSL(2,\R)$ and is invertible. Additionally, if $k(w)=\sum_{j=0}^{n-1} w_{j+1} 2^j$,  $w \in \{0,1\}^n$, then for any $x \in B_{0^n}^{n-1}$,
\begin{equation}
  \label{Lw}  L^k(x) [v]={ D \Psi_w(\Lambda_{n,F}^{-1}(x))  \Lambda_{n,F}^{-1} v \over \| D \Psi_w(\Lambda_{n,F}^{-1}(x))  \Lambda_{n,F}^{-1} v \| } \  \mod \ I,
\end{equation}
where $I$ is the reflection $I(v)=-v$.

{Recall the notation $\Omega = \{0,1\}^{\N}$} and consider the following random dynamical system $\phi$
\begin{equation}
\phi: \Z \times \R P^1 \times \Omega  \mapsto \R P^1, \quad \phi(n,[v],w):=L^n(\tau_w) [v]
\end{equation}
over a metric dynamical system $\left(\Omega, \FF, P, (p^n)_{n \in \Z}
\right)$, where $\FF$ is the cylinder
$\sigma$-algebra, generated by the cylinder sets $\{(w_1, w_2, \ldots) \in
\Omega: w_1 = a_1, \ldots, w_n = a_n \}$, $P$ is the Bernoulli measure, $p$ is
the adding machine and
\begin{equation}
  \label{tauw} \tau_{w} = \Psi^{F}_{w}(\tau)
\end{equation}
where $\tau$ is the tip of  $F$.

The skew-product $(\phi,p)$ will be denoted $\Phi$:
\[
\Phi([v],w)=(\phi(1,[v],w),p(w)).
\]
In particular, for any function $f: \R P^1 \times \Omega \mapsto \R^n$,
\[
  f \circ \Phi([v],w)= f(L(\tau_w) [v],p(w)).
\]
By Krylov-Bogolyubov theorem for a continuous dynamical system on a compact
metric space, $\Phi : \R P^1 \times \Omega \mapsto \R P^1 \times \Omega$ admits
an ergodic invariant measure $\eta^*$ (see, e.g., Section 1 of ~\cite{Sin}).

We can now apply the Birkhoff Ergodic Theorem to the following function
\begin{equation}
  \label{f_func}f([v],w) = \ln { \| D F(\tau_w) v \| \over \| v \|}.
\end{equation}
For $\eta^*$-a.a. $([v],w) \in \R P^1 \times \Omega$,
\begin{equation}
  \nonumber \lim_{N \rightarrow \infty} {1 \over N} \sum_{i=0}^{N-1} f(\Phi^i([v],w)) = \int_{\R P^1 \times \Omega} f d \eta^* =  \int_{\R P^1 \times \Omega} \ln   {  \| D F(\tau_\rho) u \| \over \| u \|   }    d \eta^*.
\end{equation}
By uniqueness of factorization for Polish spaces (see,e.g. Proposition $1.4.3$
in~\cite{Arn}), $d \eta^*=d \eta^*_w d P$, therefore, for $\eta^*$-a.a. $([v],w)
\in \R P^1 \times \Omega$,
\begin{equation}
  \label{BET} \lim_{N \rightarrow \infty} {1 \over N} \sum_{i=0}^{N-1} f(\Phi^n([v],w)) =    \int_{\Omega}  \int_{\R P^1}   \ln   {  \| D F(\tau_\rho) u \| \over \| u \|   }    d \eta_\rho^*  d P  =0
\end{equation}
by the fact that the Lyapunov exponent is equal to zero on the invariant Cantor
set.

Let $M \subset \R P^1 \times \Omega$ be a set of full measure
in $ \R P^1 \times \Omega$ for which \eqref{BET} holds, and set $M_w = M \cap
\R P^1 \times \{w\}$. By the canonicity of the Oseledets filtration the
$\Phi$-invariant bundle
\[
  E=\cup_{w \in \pi(M)} E_w, \quad E_w:={\rm span} \{M_w\},
\]
defined $\eta^*$-a.e., coincides fiber-wise with $\R P^1$: $E_w=\R P^1$.

The conclusion above about the a.e.\ zero average of the function \eqref{f_func}
has the following
consequence. 

\begin{proposition}\label{zero_average_cone_expansion}{ (\underline{Zero cone expansion.})}
  Let $k(w)=\sum_{j=0}^{\infty} w_{j+1} 2^j$, $|w|<\infty$, and $\rho \in
  \Omega$. Consider the orbits $\{v_{p^k(\rho)} = D F_*^k(\tau_\rho)
  v_{\rho}\}_{k=0}^\infty$ and $\{u_{p^k(\rho)} = D F_*^k(\tau_\rho)
  u_{\rho}\}_{k=0}^\infty$ of two non-zero vectors by $D F_*$. Denote
  $(u_w,v_w)$ the angles of these vectors. Then for a.a.
  $\rho$, $v_\rho$ and $u_\rho$
  \[
    \lim_{k \rightarrow \infty} {1 \over k} \ln { \sin(u_{p^k(\rho)}, v_{p^k(\rho)}) \over \sin(u_\rho, v_\rho)}=0.
  \]
\end{proposition}
\begin{proof}
  Consider the action of $D F_*(\tau_\rho)$ on two vectors $u_{ \rho }$ and
  $v_{\rho }$. By the fact that $D F_*( \tau_\rho )$ is a symplectomorphism,
  \begin{equation}
    \nonumber  \sin \left(u_{p^k(\rho)}, v_{p^k(\rho)}\right) =  { \| u_{ \rho  }\| \|v_{ \rho    }\| \over \|D F_*^k(  \tau_\rho ) u_{  \rho    } \|   \|D F_*^k(  \tau_\rho ) v_{ \rho    } \|  } \sin (u_{ \rho     }, v_{\rho}).
  \end{equation}
  Therefore,
  \begin{align}
    \nonumber {1 \over k} \ln { \sin \left(u_{p^k(\rho)}, v_{p^k(\rho)} \right) \over \sin \left(u_\rho, v_\rho \right)} & =  -  {1 \over k} \hspace{-0.5mm}  \ln  \hspace{-0.5mm} {\| u_{p^k(\rho)}\| \|v_{p^k(\rho)}  \| \over \|u_{\rho} \|   \|v_{\rho} \|   } \\
    \nonumber & =  - {1 \over k}  \hspace{-0.5mm} \sum_{i=0}^{k-1}  \hspace{-0.5mm} \ln  \hspace{-0.5mm} {  \| u_{p^{i+1}(\rho)}   \| \|v_{p^{i+1}(\rho)}   \|  \over  \|u_{p^i(\rho)} \|   \|v_{p^i(\rho)} \|   } \\
    \nonumber  & =   - {1 \over k}  \hspace{-0.5mm} \sum_{i=0}^{k-1}  \hspace{-0.5mm} \ln  \hspace{-0.5mm} {  \hspace{-0.3mm} \| \hspace{-0.3mm} D F_* \hspace{-0.3mm}( \tau_{p^i(\rho)}) u_{p^i(\rho)} \|  \over \|u_{p^i(\rho)} \| } \hspace{-0.5mm} -  \hspace{-0.5mm} \ln  \hspace{-0.5mm} {  \hspace{-0.3mm} \|  \hspace{-0.5mm} D F_* \hspace{-0.3mm}(  \tau_{p^i(\rho)}) v_{p^i(\rho)} \| \over \|v_{p^i(\rho)} \|   }.
  \end{align}
  The claim follows by \eqref{BET}.
\end{proof}
We will use the following notation in the next lemmas
\begin{equation}\label{Omegam}
  \Omega_m:=\left\{0^m w \in \Omega : w \in \Omega   \right\},
\end{equation}
and, we will define the following map on $\R P^{1} \times \Omega$: 
\[
  \varLambda_*([v],w)=\left( {\Lambda_*[v] \over \| \Lambda_* [v] \|} \mod I ,0w \right),
\]
where $I$ is the reflection $I(v)=-v$, and its inverse on $\R P^{1} \times {\Omega_{1}}$:
\begin{eqnarray}
 \nonumber  \varLambda_*^{-1}([v],w)&=&\left( {\Lambda_*^{-1}[v] \over \| \Lambda_*^{-1} [v] \|}  \mod I,\iota(w) \right), \\
\nonumber \iota( 0 w)&=&w.
\end{eqnarray}
Let $G:=F^2$, and consider the action of the corresponding skew-product $\Gamma$
on $\R P^1 \times \Omega_1$:
\begin{align}
  \Gamma([v],\omega)=(N(\tau_\omega)[v], p^2(\omega)), \quad   N(x) [v]={ D G(x) v \  \mod \ Z(\R^2)}.
\end{align}
It is not difficult to show that the skew-product $\Gamma_{*}
\vert_{\Omega_1}$, where $\Gamma_{*}$ denotes the map $\Gamma$
  for $G_{*}=F_{*}^{2}$, admits an ergodic invariant measure
\begin{equation}
 \label{equality} \hat \eta^*:= \varLambda_*  \eta^*= \eta_w^* \circ \Lambda_*^{-1} \times P \circ \Lambda_*^{-1}.
\end{equation}
Since $G_* \vert_{\Omega_1}$ is smoothly conjugate to $F_* \vert_{\Omega}$ (by
the linear map $\Lambda_*$), and since the Lyapunov exponents are smooth
invariants, we get, using the scaling property \eqref{equality} of $\hat \eta^*$
and the $\Phi_*$-invariance of $\eta^*$:
\begin{align*}
  0 &= \hspace{-1mm}  \int \displaylimits_{\R P^1 \times \Omega_1}  \hspace{-2mm}  \ln { \| D G_*(\tau_w) v \| \over \|v\|}  \hat \eta^*(d [v],d w) \\
  &= \hspace{-4.5mm}  \int \displaylimits_{\R P^1 \times \Omega_1} \hspace{-4.5mm}   \ln { \| D F_*(\tau_{p(w)})D F_*(\tau_w) v \|  \over  \| D F_*(\tau_w) v \|  } \hat \eta^*(d [v],d w)  + \hspace{-4.5mm}  \int \displaylimits_{\R P^1 \times \Omega_1} \hspace{-4.5mm}   \ln { \| D F_*(\tau_w) v \| \over \|v\|  } \hat \eta^*(d [v],d w) \\
  &= \hspace{-6.0mm} \int \displaylimits_{\R P^1 \times p(\Omega_1)} \hspace{-6.0mm} \ln { \| D F_*(\tau_w) v \| \over \| v\|  }  \hat \eta^* \left(\Phi_*^{-1} (d [v],d w)) \right)  + \hspace{-4.5mm} \int \displaylimits_{\R P^1 \times \Omega_1} \hspace{-4.5 mm} \ln { \| D F_*(\tau_w) v \|  \over \|v\| } \Lambda_*\hspace{-1.2mm}\left( \eta^* \right) \hspace{-0.5mm}(d [v],d w) \\
  &= \hspace{-6.0mm} \int \displaylimits_{\R P^1 \times p(\Omega_1)} \hspace{-6.0mm} \ln { \| D F_*(\tau_{0w}) \Lambda_*v \|   \over  \| \Lambda_*v \|}\Phi_*(\eta^*) (d [v],d w))  + \hspace{-5mm} \int \displaylimits_{\R P^1 \times \Omega_1} \hspace{-5.0mm} \ln { \| D F_*(\tau_{ 0 w} ) \Lambda_* v \| \over \|   \Lambda_* v \|   }  \eta^*(d [v],d w)  \\
  &= \hspace{-3mm} \int \displaylimits_{\R P^1 \times p(\Omega_1)} \hspace{-3mm} \ln { \| D \psi _1(\tau_w) v \|   \over  \| D \psi _0(\tau_w) v \|}\eta^* (d [v],d w))  + \hspace{-2mm} \int \displaylimits_{\R P^1 \times \Omega_1} \hspace{-2mm} \ln { \| D \psi _1(\tau_{ w} ) v \| \over \|   D \psi _0(\tau_w)  v \|   }  \eta^*(d [v],d w)  \\
  &= \hspace{-1mm} \int \displaylimits_{\R P^1 \times \Omega} \hspace{-2mm} \ln { \| D \psi _1(\tau_w) v \|   \over  \| v \|}\eta^* (d [v],d w))  - \hspace{-1mm} \int \displaylimits_{\R P^1 \times \Omega} \hspace{-2mm} \ln { \| D \psi _0(\tau_{ w} ) v \| \over \|    v \|   }  \eta^*(d [v],d w).
\end{align*}
We conclude that both microscope branches have the same $\eta^*$-average. We will denote this common $\eta^*$-average as  $\ln \nu$: 
\begin{align}
  \nonumber \ln \nu &=    \int \displaylimits_{\R P^1 \times \Omega}     \ln { \| D \psi _1(\tau_w) v \| \over \| v\|  }  \eta^* (d [v],d w))  \\
  \label{nu} &=  \hspace{-2mm}   \int \displaylimits_{\R P^1 \times \Omega}   \ln { \| D \psi _0(\tau_w) v \|  \over \| v \| }  \eta^* (d [v],d w)).
\end{align}
Furthermore, the Birkhoff Ergodic Theorem applied to either of the microscope
branches, $v \in \{0,1\}$, gives for a. a. $u \in \R^2 \setminus \{0\}$ and a.a.
$\rho \in \Omega$:
\begin{equation*}
  \ln \nu = \lim_{n \rightarrow \infty} {1 \over 2^n} \sum_{j=0}^{2^n-1} \ln \| D \psi_v   (\tau_{p^{j} (\rho)}) L^{j}( \tau_\rho) [v] \|.
\end{equation*}

\subsection{Adequate convergence of Birkhoff averages.} In this subsection we will prove that the Birkhoff sums \eqref{epsk} converge sufficiently fast, this convergence being sufficient, or {\it adequate}, for several stronger ergodic properties which will be proved later in this Section.


\begin{definition}{\it (\underline{Adequate convergence of Birkhoff
      averages}.)} \label{fast_conv} We will say that the Birkhoff averages of a
  function $f \in L^1(\R P^1 \times \Omega,\eta^*)$ {\it converge adequately} if
  for a.a. $([u], w) \in \R P^{1} \times \Omega$ there exist an increasing
  sequence $n_{k}$ and a sequence $\eps_{k}$, both depending on $([u], w)$, such
  that
  \begin{equation}
    \label{epsk}    \left| {1 \over 2^{n_j}} \sum_{i=0}^{2^{n_j}-1}  f \left(L^i(\tau_w)[u], p^i(w) \right) - \bar f \right|
    < \eps_{k}| \bar f|
  \end{equation}
  for all $j \ge k$, where
  \begin{equation}
    \label{barf} \bar f=\int_{\R P^1 \times \Omega} f([u],w) \eta^* d([u],w),
  \end{equation}
  and $\eps_k$ satisfies
  \begin{align*}
    \sum_{k=0}^\infty \eps_k l_k & < \infty, \\
    l_k:=n_{k+1}-n_k &< n_k.
  \end{align*}
\end{definition}

In this subsection we will prove that the Birkhoff averages converge adequately.

Define the correlation functions for random sequences
\[
  X_k^f([u],\omega)\coloneqq f\left({ D F_{*}^{k}(\tau_w) u \over \| D
      F_{*}^{k}(\tau_w) u \| },p^k(w)\right)-\bar f
\]
on $\R P^1 \times \Omega$ as
\begin{align}
\nonumber   R_k^f & = \hspace{-3mm}\int \displaylimits_{\R P^1 \times \Omega}\hspace{-2mm} \left( \hspace{-0.5mm}  f \hspace{-1mm}  \left( {u \over \| u \| },w \right)-\bar f\right) \hspace{-1mm} \left( \hspace{-0.5mm}   f \hspace{-1.0mm} \left( \hspace{-0.5mm}   {D F^k_*(\tau_w)  u  \over \|  D F^k_*(\tau_w)  u \| }, p^k(w) \right)-\bar f\right) \hspace{-0.5mm} \eta^* ( d[u], d w)\\
\label{Rkf}        & = \int \displaylimits_{\R P^1 \times \Omega}\hspace{-2mm} f\left( {u \over \| u \| },w \right) f\left(  {D F^k_*(\tau_w)  u  \over \|  D F^k_*(\tau_w)  u \| }, p^k(w) \right)-\bar f^2 \eta^* ( d[u], d w).
\end{align}
The correlation functions are a classical tool to control the speed of
convergence of Birkhoff averages. Out of many results in this direction, we will
use the following theorem from~\cite{Pet}.

\begin{theorem}[Theorem $2$ in~\cite{Pet}.]\label{Petrov}
  Let $\{X_i\}_{i=1}^\infty$ be a stationary sequence of random variables with
  $\E(X_i)=0$ and $R_k=\E(X_n X_{n+k})$. Then, for any $\delta>0$, the partial
  sums of $X_i$ satisfy:
  \[
    \sum_{k=1}^n X_i=o\left(\left( n \sum_{i=1}^{n-1} |R_i| \right)^{1 \over 2}
      (\ln n)^{{3 \over 2} +\delta } \right)
  \]
  a.s.
\end{theorem}

In  Lemmas $\ref{sum_estimates}$, $\ref{avoids}$  and  Proposition $\ref{adequate_conv}$, we show that the partial sums $\sum_{k=0}^{2^n-1}|R_k|$ of the correlation functions grow as $\kappa^n$ for some $1<\kappa<2$. According to Theorem \ref{Petrov}, this is sufficient for a fast convergence of Birkhoff averages.

Lemma $\ref{sum_estimates}$ describes how the partial sums of the correlation functions change as the orbit is advanced by one level. This is a preparatory result and the bound contains measures $\eta^*(V^n_m)$ of certain sets which will be estimated in Proposition $\ref{adequate_conv}$.

Recall the definition of $\gamma=\mu_*/|\lambda_*|$, the definition $(\ref{Omegam})$ of $\Omega_m$, and let, as before, $\beta<1$ be such that $\alpha_m^{2 n} \theta_2^{n+m}<\beta^n$ for a sufficiently large $m$ and any $n \in \N$.

Below, we will use a stand in $C$ for a positive constant whose value is immaterial to the proof.

\begin{lemma} \label{sum_estimates}
There exist cones $\Theta^n_w$ of opening $C \gamma^n$ centered on some $w$-dependent directions, such that the following holds.
  
  Set
  \begin{equation}\label{Vnm}
    V^n_m= \bigcup_{w \in \Omega_m} \Theta^n_w \times \{ w\}.
  \end{equation}

Then  the partial sums $\sum_{k=0}^{2^{n+1}-1} |R^f_k|$ of the correlation functions $(\ref{Rkf})$ satisfy the following bound
  \begin{align}
    \sum_{k=0}^{2^{n+1}-1} |R^f_k|- \sum_{k=0}^{2^{n}-1} |R^f_k| \le  2^{n+1} 2^m C \eta^*(V^n_m)+ 2^n C  (\theta_2^n+\gamma^n),
  \end{align} 
where $\eta^*$ is the Krylov-Bogolubov invariant ergodic measure and $\theta_2$ is as in Proposition $\ref{esti}$.
\end{lemma}
\begin{proof}
 Consider the partial sums $\Sigma^f_{2^{n+1}}=\sum_{k=0}^{2^{n+1}-1} |R^f_k|$,
  which we split in first $2^n$ terms and the remaining $2^n$:
  \begin{align*}
    \Sigma^f_{2^{n+1}}  &  \le   \Sigma^f_{2^n}  + \sum_{k=2^n}^{2^{n+1}-1}   \left|  \hspace{1mm}  \int \displaylimits_{\R P^1 \times \Omega}   \left( f \left( {u   \over \| u \|}, w \right)  - \bar f \right)\times \right. \\
                        & \hspace{25mm} \left. \times \left( f\left( {D F^k_*(\tau_w)  u \over \|  D F^k_*(\tau_w)  u \|   }, p^k(w)\right) -\bar f \right) \eta^* ( d[u], d w)  \right| \\
                        & \le   \Sigma^f_{2^n} +   \sum_{k=0}^{2^{n}-1}   \left|  \hspace{1mm} \int \displaylimits_{\R P^1 \times \Omega}  \left( f\left( { u \over  \| u \|   },w\right)-\bar f\right)   \times \right. \\
                        &\hspace{25mm} \left. \times \left(  f \left(  { D F^{k+2^n}_*(\tau_w)  u  \over  \|  D F^{k+2^n}_*(\tau_w)  u \|   }, p^{k+2^n}(w) \right) - \bar f \right)  \eta^* ( d[u], d w)  \right|\\
                        &=:   \Sigma^f_{2^n} + \tilde{\Sigma}^f_{2^n}.
  \end{align*}
  By the $\Phi$-invariance of the measure,
  \begin{align*}
    \tilde \Sigma^f_{2^{n}}&= \sum_{k=0}^{2^{n}-1}    \left|  \hspace{1mm} \int \displaylimits_{\R P^1 \times \Omega} \left( f \left( {D F^{-2^n}(\tau_w) u  \over  \|  D F^{-2^n}(\tau_w) u \|}, p^{-2^n}(w) \right) - \bar f \right) \times \right. \\
                           & \hspace{25mm} \left. \times \left( f \left( { D F^{k}_*(\tau_w)  u   \over \|  D F^{k}_*(\tau_w)  u \|   }, p^k(w)\right) -\bar f \right)  \eta^* ( d[u], d w) \right|.
  \end{align*}
  Again, by the invariance of the measure for any $f \in L^1(\R P^1 \times  \Omega)$ and any $m \in \N$,
  \begin{align*}
    \int \displaylimits_{\R P^1 \times \Omega} f([u],w) \eta^*(d [u],d w) &=  \sum \displaylimits_{l=0}^{2^m-1} \int \displaylimits_{\R P^1 \times p^l(\Omega_m)} f([u],w) \eta^*(d [u],d w) \\
    &= \sum_{l=0}^{2^m-1} \int \displaylimits_{\R P^1 \times \Omega_m} f(\Phi^l([u],w)) \eta^*(d [u],d w) \\
    &=\sum_{l=0}^{2^m-1} \int \displaylimits_{\R P^1 \times \Omega_m} f([u],w) \Phi^l \eta^*(d [u],d w)\\
    & =2^m \int \displaylimits_{\R P^1 \times \Omega_m} f([u],w) \eta^*(d [u],d w).
  \end{align*}
  Therefore, for any $m>n$,
  \begin{align*}
    \tilde \Sigma^f_{2^{n}} &= 2^m \sum_{k=0}^{2^{n}-1} \left|  \hspace{1mm}  \int \displaylimits_{\R P^1 \times \Omega_m} \left( f \left( { D F^{-2^n}(\tau_w) u  \over  \|  D F^{-2^n}(\tau_w) u \|   }, p^{-2^n}(w) \right) - \bar f \right) \times \right. \\
    & \hspace{25mm} \left. \times  \left( f \left({ D F^{k}_*(\tau_w)  u \over  \|  D F^{k}_*(\tau_w)  u \|   }, p^k(w) \right) - \bar f \right)  \eta^* ( d[u], d w) \right|.
  \end{align*}
  By the fact that $F_*$ is a renormalization fixed point, and by reversibility
  of $F_*$,
  \begin{align}
    \nonumber   \tilde \Sigma^f_{2^{n}} &= 2^m\sum_{k=0}^{2^{n}-1}   \left|  \hspace{1mm}  \int \displaylimits_{\R P^1 \times \Omega_m} \hspace{-3mm} \left( \hspace{-0.5mm}  f  \hspace{-0.5mm}\left( { T \Lambda^n_* D F_*(T(\tau_{\iota^n(w)})) \Lambda^{-n}_* T u  \over   \|   T \Lambda^n_* D F_*(T(\tau_{\iota^n(w)})) \Lambda^{-n}_*  T u \|   },  p^{-2^n}(w)\right) - \bar f \right)  \times \right. \\
    \label{int2} & \hspace{30mm} \left. \times \left( f \left( { D F^{k}_*(\tau_w)  u  \over  \|  D F^{k}_*(\tau_w)  u \|   }, p^k(w) \right) - \bar f \right)  \eta^* ( d[u], d w) \right|.
  \end{align}
Notice, that by the ratchet phenomenon (see Subsection $\ref{subsectionratchet}$), $[DF_*(T(\tau_{\iota^n(w)}))]^{-1}$ maps the vertical vector $h$ into some horizontal cone, therefore the angle between the vector  $[DF_*(T(\tau_{\iota^n(w)}))]^{-1} h$ and the vertical direction is bounded away from $0$ independently of $n$. Furthermore, $T \Lambda_*^n$ maps  $[DF_*(T(\tau_{\iota^n(w)}))]^{-1} h$ into a horizontal cone of some opening $C \gamma^n$. Therefore, there exists $N \in \N$ such that for any $n>N$ there exists a  cone
  \begin{equation*}
    \Theta^n_w:=\Theta_{ T \Lambda_*^n [DF_*(T(\tau_{\iota^n(w)}))]^{-1} h}(C \gamma^n),
  \end{equation*}
  and an integer $m>n$, such that for all $w \in \Omega_m$  and $[u] \in \R P^1 \setminus \Theta^n_w$,
  \[
    [T \Lambda^n D F _* ( T(\tau_{\iota^n(w)}) ) \Lambda^{-n} _* T u] \in
    \Theta_e(C \gamma^n),
  \]
  where, as before, $e=[1,0]$, $h=[0,1]$ and $\gamma=\mu_*/|\lambda_*|$. Indeed, 
  the cone $\Theta^n_w$ is mapped by $D F _* (T(\tau_{\iota^n(w)}))
  \Lambda^{-n}_* T$ into a cone which contains some cone $\Theta_h(\alpha_0)$
  centered on the vertical unit vector $h$, where $\alpha_0$ is bounded away from $0$ independently of $n$. Then $T \Lambda_*^n$ maps  $\Theta_h(\alpha_0)$ into a cone whose complement is $\Theta_e(C \gamma^n)$.

  We have, for any  $f \in C(\R P^1 \times \Omega)$ such that $f([e],0^\infty) \ne 0$, any unit $u \in \Theta_e(C \gamma^n)$ and any $\rho \in \Omega$:
  \begin{equation*}
   (1  - C_f \theta_2^n) (1 - C_f \gamma^n)  <  {|f([u],0^n \rho)|\over  |f([e],0^\infty)| } < (1 + C_f \theta_2^n) (1+ C_f \gamma^n),
  \end{equation*}
for some $C_f$, dependent on $f$, i.e. the first factor in \eqref{int2} is constant up to a  geometrically small quantity.

Recall the definition  $(\ref{Vnm})$  of $V^n_m$, and let $U^n_m=\left( \R P^1 \times \Omega_m \setminus V^n_m \right)$ be its complement in $\R P^1 \times \Omega_m$. We have
  \begin{align*}
    \tilde \Sigma^f_{2^{n}}  &=  2^m  \sum_{k=0}^{2^{n}-1}    \left|  \hspace{1mm}  \int \displaylimits_{U^n_m}  \left( C + O(\theta_2^n + \gamma^n) \right)  \times \right. \\
                             & \left. \hspace{20mm}  \times \left( f \left({ D F^{k}_*(\tau_w)  u  \over \|  D F^{k}_*(\tau_w)  u \|   }, p^k(w)   \right) - \bar f \right)  \eta^* ( d[u], d w) \right|+ \\
                             & \hspace{20mm} +C  2^m\sum_{k=0}^{2^{n}-1} \eta^*\left( V^n_m  \right) \\
                             & \le  2^m \sum_{k=0}^{2^{n}-1}  C  \left|  \hspace{1mm}  \int \displaylimits_{U^n_m}   \left( f \left( {  D F^{k}_*(\tau_w)  u  \over \|  D F^{k}_*(\tau_w)  u \|   }, p^k(w) \right) - \bar f \right) \eta^* ( d[u], d w) \right|+\\
                             & \hspace{20mm}+ 2^n 2^m  O(\theta_2^n+\gamma^n) +  C 2^n 2^m \eta^*\left( V^n_m \right).
  \end{align*}
  The deviation of the first integral from $0$ is of order $\eta^*(V^n_m)$.
  Therefore,
  \begin{equation}
    \label{tildeS} \tilde \Sigma^f_{2^{n}} \le  2^{n+m} \left( C \eta^*(V^n_m)+ C  (\theta_2^n+\gamma^n)+ C  \eta^*\left(V^n_m \right) \right).
  \end{equation}
\end{proof}

In the next Proposition we will estimate $\eta^*(V^n_m)$ and complete the proof of the fast convergence of Birkhoff averages. 

\begin{proposition}(\underline{Adequate convergence.})  \label{adequate_conv}    The Birkhoff averages of a function $f \in C(\R P^1 \times \Omega)$ converge   adequately with $l_k=1$, $\{n_k\}_{k=1}^\infty=\{0\} \cup \N$ and
  \[
    \eps_k \equiv \epsilon_n < C \sigma^n n^{{3 \over 2}+\delta}
  \]
  for some $0<\sigma <1$ and any $\delta>0$.
\end{proposition}

\begin{proof}
  Notice that by convergence of Birkhoff averages, for any fixed natural $n$ and
  $m$, $m>n$, a.a. unit  $v \in \R^2$, and any $\delta>0$ there exists $N=N(\delta,n,m,v) \in \N$, such
  that for any $k>N$,
  \begin{align}
    \label{etaBirk}  \eta^*\left(V^n_m \right) & < \delta+ {1 \over 2^k} \sum_{w \in \{0,1\}^{k}} \chi_{V^n_m} \left( { D \Psi_w (\tau) \Lambda_*^{-k} v  \over  \|   D \Psi_w (\tau) \Lambda_*^{-k} v   \|  }, w  \right) \\
    \nonumber  & = \delta+ {1 \over 2^k} \sum_{\rho \in \{0,1\}^{k-m}} \chi_{\Theta^n_{0^m \rho}} \left( { D \Psi_{0^m \rho} (\tau) \Lambda_*^{-k} v  \over  \|   D \Psi_{0^m \rho} (\tau) \Lambda_*^{-k} v    \|  } \right),
  \end{align}
  where $\chi_A$ is the characteristic function of the set $A$.
  
  Fix $m=2 n$ and a unit vector $v \in \R^2$ for which the Birkhoff sum in $(\ref{etaBirk})$ converges. Also, fix a number $0<\beta_1<1$
  and $\delta=2^{-2 n} \beta_1^n$. Then any vector which lies outside a vertical
  cone $\Theta_h(C \alpha^n)$ of opening $C \alpha^n$ for any fixed $\gamma<\alpha<1$ is mapped by $D
  \Psi_{0^m}=\Lambda_*^{2 n}$ outside of  $\Theta^n_{0^m \rho}$. Indeed, the image of any
  vector $(\alpha^n,1)$ under $\Lambda_*^{2 n}$ is $(\lambda_*^{2 n} \alpha^n,
  \mu_*^{2 n})$, which lies in the cone $\Theta_e\left(C \gamma^{2 n}/\alpha^n
  \right)$ whose opening is much less than $C \gamma^n$ if $\gamma/\alpha<1$. This
  cone does not intersect $\Theta^n_{0^m \rho}$.

  We will ask the following question now: How often is the vector
  $\Lambda_*^{-k} v$ mapped to the cone $\Theta_h(C \alpha^n)$ by $D
  \Psi_\rho(\tau)$ with $\rho \in \{0,1\}^{k-2n}$?

  Fix $v=e$, the horizontal unit vector (if the Birkhoff sum does not converge for this choice of $v$, then take an arbitrarily small perturbation of $e$). It is not too difficult to identify a subset of $\rho$'s for which the horizontal vector $\Lambda_*^{-k} e$ is {\it not} mapped to the cone $\Theta_h(C \alpha^n)$ by $D  \Psi_\rho(\tau)$. We will require the following short Lemma:

  \begin{lemma} \label{avoids}
    There exists a neighborhood $\UU$ of the tip $\tau$, a horizontal cone $\Theta_e$ and a vertical cone $\Theta_h$, such, that for any $x \in \UU$ and unit $v \in \Theta_e$,  $D F_*(x) v$ avoids $\Theta_h$.
  \end{lemma}
  \begin{proof}
    According to $(\ref{sdef})$, the action of the derivative at the tip on a horizontal vector is given by
    \begin{align*}
      D F_*(\tau) e&=\left[
      \begin{array}{c c}
        -{s_2^*(y(0,0),0) \over s_1^*(y(0,0),0)} & -{1 \over s_1^*(y(0,0),0)} \\
        (s_1^*(0,y(0,0)-s_2^*(0,y(0,0)){s_2^*(y(0,0),0) \over s_1^*(y(0,0),0)} & -{s_2^*(0,y(0,0)) \over s_1^*(y(0,0),0)}
      \end{array}
      \right]
      \left[
      \begin{array}{c}
        1 \\
        0
      \end{array}
      \right]
      \\
      &=
      \left[
        \begin{array}{c}
        -{s_2^*(y(0,0),0) \over s_1^*(y(0,0),0)}  \\
        (s_1^*(0,y(0,0)-s_2(0,y(0,0)){s_2^*(y(0,0),0) \over s_1^*(y(0,0),0)}.
      \end{array}
        \right]
    \end{align*}
    Since $y(0,0)=1$ and $s_1^*(1,0)=s_1(0,1)^*=1$ (recall the symmetry condition for $s^*$ from Theorem $\ref{hyp}$), we get that $\pi_1 D F_*(\tau) e=-s_2^*(1,0)$. We can now differentiate the fixed point equation $s^*(x,y)=\mu_*^{-1} s^*(z^*(x,y),\lambda_* y)$  (c.f.  $(\ref{Rs})$) 
   with respect to $y$, to get
  $$s_2^*(x,y)=\mu^{-1}_* s_1^*(z^*(x,y),\lambda_* y) z_2^*(x,y) +\mu^{-1}_* \lambda_* s_2^*(z^*(x,y),\lambda_* y),$$
  and evaluate this expression at $(x,y)=(1,0)$:
  $$s_2^*(1,0)=\mu^{-1}_* s_1^*(1,0) z_2^*(1,0) +\mu^{-1}_* \lambda_* s_2^*(1,0)=\mu^{-1}_* z_2^*(1,0) +\mu^{-1}_* \lambda_* s_2^*(1,0).$$
We can see, that $s_2^*(1,0)$ is non-zero whenever $z_2^*(1,0)=z_1^*(0,1)$ is non-zero. To demonstrate that $z_2^*(1,0) \ne 0$, we write the fixed point equation, now evaluated at $(y,x)$, once again:
$$s^*(y,x)=\mu^{-1}_* s^*(z^*(y,x),\lambda_* x),$$
differentiate it with respect to $y$:
$$s^*_1(y,x)=\mu^{-1}_* s^*_1(z^*(y,x),\lambda_* x) z_1^*(y,x),$$
and evaluate the last expression at point $(y,x)=(0,1)$, using the symmetry condition for $s^*$: $s_1^*(0,1)=s_1^*(1,0)=1$,
$$1=\mu^{-1}_* s^*_1(1,\lambda_*) z_1^*(0,1).$$
The last identity demonstrates that $z_1^*(0,1) \ne 0$.

We have obtained that $ D F_*(\tau) e$ has a non-zero horizontal component, and so does $ D F_*(x) v$ for any other unit vector $v$ in a sufficiently narrow horizontal cone and an{k-2 n-s+1}y $x$ in a sufficiently small neighborhood of $\tau$.
  \end{proof}

 Thus, $D F_*(x) v$, with $v$ in a certain horizontal cone, avoids the  vertical cone $\Theta_h(C \alpha^n)$ for all points in the Cantor set in some neighborhood of the tip if $n$ is sufficiently large. Choose an integer $s$, and  $k$ such that $k-m=k-2 n$ (recall that $m=2n$) is divisible by $s$. Consider the orbit of a horizontal vector $\Lambda_*^{-k} e$ by $D  \Psi_{\rho_i}(\tau_{\rho_{i+1} \ldots \rho_{k -2 n}})$ with $\rho \in \{0,1\}^{k-2 n}$ given by 
 \begin{equation*}
   \rho=(\rho_1, 0, \ldots, 0, \rho_{s+1}, 0, \ldots 0, \rho_{2 s+1}, 0, \ldots, 0,  \rho_{k-2 n-s+1}, 0, \ldots, 0).
 \end{equation*}
  If $s$ is sufficiently large, then the pieces $\Lambda^{s-1}_*$ of the microscope of length $s-1$ map a vector which avoids a definite vertical cone into a definite horizontal cone. If $s$ is sufficiently large, this horizontal cone has a narrow opening, and, a further application of  $D  \Psi_{\rho_{i s+1}}(\tau_{0 \ldots \rho_{k -2 n}}))$, which is equal either to $\Lambda_*$ or $D F_*(\tau_{00 \ldots \rho_{k -2 n}}) \Lambda_*$, $i=0,1, \ldots$, according to Lemma $\ref{avoids}$, again maps it outside of a definite vertical cone (if $s$ is chosen large enough, so that $\tau_{0 \ldots \rho_{k -2 n}} \in \UU$, the neighborhood from Lemma $\ref{avoids}$). We conclude that the orbit of   $\Lambda_*^{-k} e$ for such $\rho$'s never enters the vertical cone $\Theta_h(C \alpha^n)$. The number of such $\rho$'s is $2^{{k-2 n} \over s}$

  As the result of the above discussion, we conclude that there exists $s$, such that the number of words  $\rho \in \{0,1\}^{k-2 n}$ for which   $D \Psi_{\rho}(\tau) \Lambda_*^{-k} e \cap \Theta_h(C \alpha^n) \ne \emptyset$ is at most $2^{k-2 n -{{k-2 n} \over s}}=2^{{s -1 \over s}(k-2 n)}$. Therefore,
    \begin{equation*}
    \eta^*\left(V^n_{2 n} \right) \le 2^{-2 n} \beta_1^n + {1 \over 2^k} 2^{{s-1 \over s} (k -2 n)}  < 2^{-2 n} \left( \beta_1^n + 2^{{2 \over s} n-{k \over s} } \right).
  \end{equation*}
    Recall, that $k$ has to be chosen larger than  $N=N(2^{-2 n} \beta_1^n,n,2 n,e)$ for  $(\ref{etaBirk})$ to hold. In particular,  it can be chosen large enough so that $2/s-{k/( n s)}<-1$. Then
    \begin{equation*}
    \eta^*\left(V^n_{2 n} \right) < 2^{-2 n} \left( \beta_1^n + 2^{-n} \right).
  \end{equation*}
    Choose $\beta_2= \max\{\beta_1, 2^{-1}\}$, then
  \begin{equation*}
    \eta^*\left(V^n_{2 n} \right) \le C 2^{-2 n} \beta_2^n.
  \end{equation*}
  Plugging this into $(\ref{tildeS})$ with $m=2 n$, and using the fact that, according to $(\ref{theta2})$ and $(\ref{lambdamustar})$, $(\theta_2^n+\gamma^n) < C 2^{-2 n} \beta_3^n$ for some $\beta_3<1$,  we obtain that there exists $1<\kappa<2$,
  such that
  \begin{equation*}
    \tilde \Sigma^v_{2^{n}} \le C \kappa^n \implies \Sigma_{2^n}^v \le { C \over \kappa-1}  \kappa^{n+1}.
  \end{equation*}

  By Theorem \ref{Petrov}, for almost all $([u],\rho) \in \R P^1 \times \Omega$,
  \begin{equation} \label{epsbound}
    \left| {1 \over 2^{n}} \sum_{k=0}^{2^n-1} f \left( { D F_*^k(\tau_\rho) u  \over  \|   D {F_*^k}(\tau_\rho)  u    \|  }, p^k(\rho) \right) - \bar f \right| \le g([u], \rho) \  \left({\kappa \over 2} \right)^{n/2} n^{{3 \over 2}+\delta},
  \end{equation}
  where $g$ is some function on $\R P^1 \times \Omega$. Therefore, \eqref{epsk}
  holds with $l_k=1$, $\{n_k\}_{k=1}^\infty=\N$, and a summable sequence
  \begin{equation*}
    \epsilon_n < C  \left({\kappa \over 2} \right)^{n/2} n^{{3 \over 2}+\delta}.
  \end{equation*}
  The claim follows.
\end{proof}

\subsection{Bounded distortion of a microscope branches in a piece.}

\comment{
We begin with a bound on the average of the angle between pairs of vectors in the orbits by the derivative cocycle. 

\begin{lemma}\label{averageangle}
  For a.a. $\rho \in \Omega_n$  and a.a  non-zero vectors $u$ and $v$, the average log-expansion of angles between the vectors in the orbit of length $2^{n}$, 
  \begin{equation}
   \nonumber  \mathcal{A}_{n}^\rho={1 \over 2^{n}} \sum_{k=0}^{2^{n}-1}  \ln {\sin \left(DF^k_*(\tau_{\rho}) u, \ DF^k_*(\tau_{\rho}) v \right)   \over \sin (u, v)  },
  \end{equation}
where $\tau_\rho \in B_{0^n}^{n-1}$,  satisfies
  \begin{equation}
    \nonumber \mathcal{A}_{n}^\rho \le  {  \sin\left( \Lambda_*^{-{n}} u,  \Lambda_*^{-{n}} v  \right) \over  \sin(u,v)  } \left(n \ln |\lambda_*| \mu - 2 \left(n +\sum_{j=0}^{n-1} \epsilon_j \right) \ln \nu \right)
  \end{equation}
\end{lemma}
\begin{proof}
  Using the notation of \eqref{barrho} and \eqref{tilderho}, and denoting
  $$\hat{\mathcal{A}}_{n}^\rho=\mathcal{A}_{n}^\rho { \sin(u,v) \over \sin\left( \Lambda_*^{-{n}} u,  \Lambda_*^{-{n}} v  \right)},$$
  we get
  \begin{align}
    \nonumber   \hat{\mathcal{A}}_{n}^\rho&=  {1 \over 2^{n}}  \sum_{w \in \{0,1\}^{n}}  \ln { \sin \left(D \Psi_w (\tau_{ \iota^n(\rho) }) \Lambda_*^{-{n}} u,  D \Psi_w  (\tau_{\iota^n(\rho) }) \Lambda_*^{-{n}} v\right) \over \sin \left(\Lambda_*^{-{n}} u,  \Lambda_*^{-{n}} v\right)   } \\
    \nonumber &= n \ln |\lambda_*| \mu - {1 \over 2^{n}}  \sum_{w \in \{0,1\}^{n}} \ln { \left\|  D \Psi_w (\tau_{\iota^n(\rho) }) \Lambda_*^{-{n}} u \right\| \over \|  \Lambda_*^{-{n}}  u\|} { \left\|  D \Psi_w  (\tau_{\iota^n(\rho) }) \Lambda_*^{-{n}} v \right\| \over  \| \Lambda_*^{-{n}}  v \|   }    \\
    \nonumber &= n \ln |\lambda_*| \mu - {1 \over 2^{n}}  \sum_{w \in \{0,1\}^{n}} \ln \prod_{i=1}^{{n}} { \left\|  D \Psi_{\tilde w_i} (\tau_{\iota^n(\rho) }) \Lambda_*^{-{n}} u \right\| \over \left\| D \Psi_{\tilde w_{i+1}} (\tau_{\iota^n(\rho) }) \Lambda_*^{-{n}} u \right\| } + \\
    \nonumber & \hspace{20.3mm} - {1 \over 2^{n}}  \sum_{w \in \{0,1\}^{n}} \ln \prod_{i=1}^{{n}} { \left\|  D \Psi_{\tilde w_i} (\tau_{\iota^n(\rho) }) \Lambda_*^{-{n}} v \right\| \over \left\| D \Psi_{\tilde w_{i+1}} (\tau_{\iota^n(\rho) }) \Lambda_*^{-n} v \right\| }  \\
    \nonumber & =n \ln |\lambda_*| \mu + \Sigma_u +\Sigma_v,
  \end{align}
  where
  \begin{align}
    \nonumber \Sigma_u &= - {1 \over 2^{n}}  \sum_{i=1}^{n} \sum_{\bar w_{i-1} \in \{0,1\}^{i-1} \atop \tilde w_{i+1} \in \{0,1\}^{{n}-i} } \ln { \left\|  D \psi_0(\tau_{\tilde w_{i+1} \iota^n(\rho) })  D \Psi_{\tilde w_{i+1}} (\tau_{\iota^n(\rho) }) \Lambda_*^{-{n}} u \right\| \over \left\| D \Psi_{\tilde w_{i+1}} (\tau_{\iota^n(\rho) }) \Lambda_*^{-{n}} u \right\| } + \\
    \nonumber & \hspace{35mm}+ \ln { \left\|  D \psi_1(\tau_{\tilde w_{i+1} \iota^n(\rho) })  D \Psi_{\tilde w_{i+1}} (\tau_{\iota^n(\rho) }) \Lambda_*^{-{n}} u \right\| \over \left\| D \Psi_{\tilde w_{i+1}} (\tau_{\iota^n(\rho) }) \Lambda_*^{-{n}} u \right\| } \\
    \nonumber &\le  - {1 \over 2} \sum_{i=1}^{n} {1 \over 2^{i-1}}\sum_{\bar w_{i-1} \in \{0,1\}^{i-1}} (1 +\epsilon_{{n}-i}) \ln \nu +  (1 +\epsilon_{{n}-i}) \ln \nu \\
    \nonumber &\le  - \ln \nu \sum_{i=1}^{n}(1 +\epsilon_{{n}-i}) \\ 
    \nonumber &=  - \left(n +\sum_{j=0}^{n-1} \epsilon_{j} \right) \ln \nu. 
  \end{align}
  and similarly for $\Sigma_v$.

  Notice that the sequence $\epsilon_j$ is summable. The claim follows.
\end{proof}
}

We will now study the  bounds on the average ratio of actions of derivatives of long iterates on a vector, evaluated at two different points in a deep piece. 

 Let $\eta$ be a non-zero unit vector. Set
\begin{align}
  \label{u_vecs}  u_w  & = D F^k(x_{0^n}) \eta= D \Psi_w (x) \Lambda_*^{-n} \eta, \\
  \label{v_vecs}  v_w &= D F_*^k (y_{0^n}) \eta= D \Psi_w (y)  \Lambda_*^{-n} \eta,
\end{align}
where $x= \Lambda^{-n}_* x_{0^n}$. Define the average log-distortion of the derivative cocycle as
\begin{align}
  \nonumber   \delta_{n}(x_{0^n},y_{0^n},\eta) & =  {1 \over 2^{n}} \sum_{k=0}^{2^n-1} \ln { \|  D F^k_*(x_{0^n}) \eta\| \over \|  D F^k_*(y_{0^n}) \eta\| } \\
  \label{deltal} &=    {1 \over 2^{2 n}} \sum_{w \in \{0,1\}^n}  \ln { \| u_w \|  \over \| v_w \|  }.
\end{align}

Our last result is that the log-distortion is bounded.

\begin{proposition}\label{angles3} { (\underline{Distortion of the derivative cocycle
      in a piece.})} There exist constants $C>0$ and $0<c<1$ such that the   following holds for any $n \in \N$, any two $x_{0^n}, y_{0^n} \in   B_{0^{n}}^{n-1}$ and a.a non-zero vectors $\eta$, the log-distortion $\delta_n\equiv \delta_{n}(x_{0^n},y_{0^n},\eta)$ satisfies
\begin{equation*}    \ln \hspace{-0.0mm} \left\{ {  1 \hspace{-0.0mm}-\hspace{-0.0mm}  c   \over   1 \hspace{-0.0mm}+\hspace{-0.0mm}  C    } \nu^C \right\} \hspace{-0.0mm}-\hspace{-0.0mm} {1 \over 2} \ln{ \left(1\hspace{-0.0mm}+\hspace{-0.0mm} C \|x \hspace{-0.0mm}-\hspace{-0.0mm} y\| \right)} \le    \delta_{n} \le \ln  \hspace{-0.0mm} \left\{ {  1\hspace{-0.0mm}+\hspace{-0.0mm}  C  \over   1\hspace{-0.0mm}-\hspace{-0.0mm}  c   } \nu^{-C} \right\} \hspace{-0.0mm}+\hspace{-0.0mm}{1 \over 2} \ln{ \left(1\hspace{-0.0mm}+\hspace{-0.0mm} C \|x\hspace{-0.0mm}-\hspace{-0.0mm}y\| \right)},    
\end{equation*}
where  $x= \Lambda^{-n}_* x_{0^n}$ and  $y= \Lambda^{-n}_* y_{0^n}$.
\end{proposition}
\begin{proof}
  Recall that \eqref{epsk}, with the function
  \begin{equation}
    f_v([\eta],\rho)=\ln { D \psi_v (\tau_\rho) \eta \over \|\eta\|}, \ v=0,1,
  \end{equation}
  is satisfied for almost all $([u],w)$. Choose a point $\tilde \tau \in  \CC_{F_*}$ so that \eqref{epsk} is satisfied for $\tilde \tau_{0^n}$ and a.a. $\eta$.

For such $\tilde \tau \in \C_{F_*}$, define the following sequence:
\begin{align}
  \label{uwi}  \tilde u_{\tilde w_i}  & =D \Psi_{\tilde w_i} (\tilde \tau) \Lambda_*^{-n} \eta,  \quad 1 \le i \le n.
\end{align}
where $\tilde w_{i}$  has been defined in $(\ref{tilderho})$. 
Similarly, define, for a non-zero vector $z$, a sequence
  \begin{equation}
    \label{zwi} z_{\tilde w_{i}}  = D  \Psi_{\tilde w_i} (\tilde \tau) z, \quad 1 \le i \le n.
  \end{equation}
The vectors $\tilde{u}_{\tilde w_{i}}$ are in the forward orbit of the vector $\Lambda_*^{-n} \eta$ under a microscope branch evaluated at a different point than the forward orbit $ D \Psi_{\tilde w_i}(x)\Lambda_*^{-n} \eta$.
  
Then, 
  \begin{equation}
    \label{sigmas}   {1 \over 2^{n}}  \hspace{-1mm}  \sum_{ { w \in \{0,1\}^{n}}}   \hspace{-3mm}   \ln {\| u_{w}  \| \over \| \Lambda_*^{-n} \eta \| } =   {1 \over 2^{n}}  \hspace{-1mm}  \sum_{ { w \in \{0,1\}^{n}}}   \hspace{-3mm}   \ln {\| \tilde u_{w}  \| \over \| \Lambda_*^{-n} \eta \| } +  \ln {\| u_w \| \over \| \tilde u_{w}\| }  =  \Sigma_1+\Sigma_2,
  \end{equation}
  where
  \begin{equation}
    \label{sigma1}  \Sigma_1 ={1 \over 2^{n}} \hspace{-1mm}  \sum_{{ w \in \{0,1\}^{n}}}   \sum_{i=1}^{n}   \ln   {\|   D \psi_{w_i} (\tilde \tau_{\tilde w_{i+1}})  D \Psi_{\tilde w_{i+1}} (\tilde \tau)\Lambda^{-n}_{*} \eta \|  \over  \|   D \Psi_{\tilde w_{i+1}}(\tilde \tau)\Lambda^{-n}_{*} \eta   \|},
  \end{equation}
with the convention,  $D \Psi_{\tilde w_{n+1}}(\tilde \tau)=\mathbb{I}$,  and
  \begin{equation}
    \nonumber \Sigma_2={1 \over 2^{n}}  \sum_{{ w \in \{0,1\}^{n }}}   \ln {\| u_w \| \over \| \tilde u_{w}\| }.
  \end{equation}
By Lemma $\ref{more_bounds_direct}$,
\begin{align}
  \nonumber \Sigma_2 & \le {1 \over 2^{n}}  \sum_{{ w \in \{0,1\}^{n}}}  \ln \Pi_{1}  \prod_{j=n}^{1}    \left( 1  + K\theta_2^{n\shortminus j} \right)^{  \hspace{-0.5mm} |w_{j}|} \\
  \nonumber & \le {1 \over 2^{n}}  \sum_{ w \in \{0,1\}^{n}}   \ln  \left( 1+  C \right) + \Sigma_2' \\
  \label{Sigma2} & \le \ln  \left( 1+  C \right) + \Sigma_2',
\end{align}
where we have used that
$$\prod_{j=n}^{1}  \hspace{-0.3mm}   \left( 1 \hspace{-0.6mm}   + \hspace{-0.6mm}  K\theta_2^{n\shortminus j} \right)^{  \hspace{-0.5mm} |w_{j}|}  \le  1 + C$$
for some $C>0$, and have defined 
\begin{equation*}
  \Sigma_2'= {1 \over 2^{n}}  \sum_{{w \in \{0,1\}^{n}}}  \ln  {\sin\left({u}_{w_n},z_{w_n}  \right) \over \sin\left(\tilde{u}_{w_n}, z_{w_n}  \right) }\  { \sin\left( \tilde{u}_w, z_{w} \right) \over \sin \left(u_w, z_{w} \right)}.
\end{equation*}
Next, we address $\Sigma_1$. The sum in  \eqref{sigma1} is of the form
  \begin{equation*}
    \Sigma_1\coloneqq{1 \over 2^{n}}  \sum_{i=1}^{n} \sum_{\bar w_{i-1} \in \{0,1\}^{i-1},  \atop \tilde w_{i+1} \in \{0,1\}^{n-i}}   g_0(\tilde w_{i+1})+g_1(\tilde w_{i+1}),
  \end{equation*}
  where $g_0(\tilde w_{i+1})$ and $g_1(\tilde w_{i+1})$ are the quantities (corresponding to  the choices $w_i=0$ and $w_i=1$) that depend on $\tilde w_{i+1}$ only but not on $\bar w_{i-1}$. Therefore,
  \begin{align*}
    \Sigma_1 &= {1 \over 2^{n}}  \sum_{i=1}^{n}  \sum_{\bar w_{i-1} \in \{0,1\}^{i-1}}   \sum_{\tilde w_{i+1} \in \{0,1\}^{n-i}}  g_0(\tilde w_{i+1})+g_1(\tilde w_{i+1}) \\
      &= \sum_{i=1}^{n} {1 \over 2^{n-i+1}}  \sum_{\tilde w_{i+1} \in \{0,1\}^{n-i}}     g_0(\tilde w_{i+1})+g_1(\tilde w_{i+1})  \\
    &= {1 \over 2} \sum_{i=1}^{n} \left( {1 \over 2^{n-i}}  \sum_{\tilde w_{i+1} \in \{0,1\}^{n-i}}     g_0(\tilde w_{i+1})+g_1(\tilde w_{i+1}) \right).
  \end{align*}  
    From here we see that we can apply \eqref{epsk} for the sums in parenthesis, twice, once for $g_{0}$ and  once for $g_{1}$, taking into account that by \eqref{epsbound}, each sum in parenthesis deviates from $\ln \nu$ by no more than $\epsilon_{n-i} |\ln \nu|$.  Therefore, since $\ln\nu < 0$,
  \begin{equation}
 \label{Sigma1}   \Sigma_1 \leq \sum_{j=0}^{n-1} (\ln\nu - \epsilon_{j} \ln\nu ) =  \left(n -\sum_{j=0}^{n-1} \epsilon_j \right)\ln \nu.
  \end{equation}

`  Putting this into \eqref{sigmas}, we get
  \begin{equation}
\label{exp1}    {1 \over 2^{n}}  \hspace{-1mm}  \sum_{ { w \in \{0,1\}^{n}}}   \hspace{-3mm}  \ln {\| u_{w}  \| \over \| \Lambda_*^{-n} \eta \| }   \le \ln \left\{ \left( 1+  C \right) \nu^{ n -\sum_{j=0}^{n-1} \epsilon_j} \right\} + \Sigma_2'.    
  \end{equation}

  We can now repeat similar calculation for the sequence $v_w$ to obtain a bound from below:
  \begin{equation}
    \label{exp2}  {1 \over 2^{n}}  \hspace{-1mm}  \sum_{ { w \in \{0,1\}^{n}}}   \hspace{-3mm}   \ln { \| v_{w}  \|   \over \| \Lambda_*^{-n} \eta \| }  \ge \ln \left\{ \left( 1-  c \right) \nu^{ n +\sum_{j=0}^{n-1} \epsilon_j} \right\} + \Sigma_2'',
  \end{equation}
  where
\begin{equation*}
  \Sigma_2''={1 \over 2^{n}}  \sum_{{w \in \{0,1\}^{n}}}  \ln  {\sin\left(v_{w_n},z_{w_n}  \right) \over \sin\left(\tilde{u}_{w_n}, z_{w_n}  \right) }\  { \sin\left( \tilde{u}_w, z_{w} \right) \over \sin \left(v_w, z_{w} \right)}.
\end{equation*}
We, therefore, obtain for $ \delta_{n}(x_{0^n},y_{0^n},\eta)$:
  \begin{equation*}
     \delta_{n}(x_{0^n},y_{0^n},e) \le \ln \left\{ {  1+  C  \over   1-  c  } \nu^{-2 \sum_{j=0}^{n-1} \epsilon_j} \right\} + \Sigma_2'-\Sigma_2'',  
  \end{equation*}
  where
  \begin{align*}
    \Sigma_2'-\Sigma_2''  &= {1 \over 2^{n}}  \sum_{{w \in \{0,1\}^{n}}}  \ln  {\sin\left({u}_{w_n},z_{w_n}  \right) \over \sin\left(v_{w_n}, z_{w_n}  \right) }\  { \sin\left( v_w, z_{w} \right) \over  \sin \left( u_w, z_{w} \right)}.
  \end{align*}
  We can now choose $z$ so that
  $$\sin\left(u_{w},z_{w}  \right)=\sin\left(v_{w},z_{w}  \right).$$
  Then
\begin{align*}
  \Sigma_2'-\Sigma_2'' & = {1 \over 2^{n}}  \sum_{{w \in \{0,1\}^{n}}}  \ln  {\sin\left({u}_{w_n},z_{w_n}  \right) \over \sin\left(v_{w_n}, z_{w_n}  \right) }\\
  & = {1 \over 2^{n}}  \sum_{{w \in \{0,1\}^{n-1}}}  \ln  {\sin\left({u}_0,z_0  \right) \over \sin\left(v_0, z_0  \right) } +  \ln  {\sin\left({u}_1,z_1  \right) \over \sin\left(v_1, z_1  \right) }\\
  & = {1 \over 2^{n}}  \sum_{{w \in \{0,1\}^{n-1}}}  \ln  {\sin \left( \Lambda_*^{-n+1} \eta ,z_0  \right) \over \sin\left( \Lambda_*^{-n+1} \eta , z_0  \right) } +  \ln  {\sin\left( D F_*(x_0)  \Lambda_*^{-n+1} \eta,z_1  \right) \over \sin\left( D F_*(y_0)  \Lambda_*^{-n+1} \eta, z_1  \right) }\\
     & = {1 \over 2} \ln{ \left(1+ C \|x-y\| \right)}.
\end{align*}
We conclude that
  \begin{align*}
    \delta_{n}(x_{0^n},y_{0^n},e) &\le \ln \left\{ {  1+  C \over   1-  c  } \nu^{-2 \sum_{j=0}^{n-1} \epsilon_j} \right\} + {1 \over 2} \ln{ \left(1+ C \|x-y\| \right)}.
  \end{align*}
  Since $\epsilon_j$ is summable, the power of $\nu$ is bounded. The claim follows.
  
  To obtain a lower bound one can follow the steps of the above proof, using the lower bound from Lemma $\ref{more_bounds_direct}$ for $\Sigma_2$ in $(\ref{Sigma2})$, changing the sign in front of $\epsilon_{j} \ln\nu$ in $(\ref{Sigma1})$, and further, in front of the constant  in the expressions
$(\ref{exp1})$ and $(\ref{exp2})$.

\end{proof}

The above proposition means that the average log-distortion is bounded as epressed by the inequality $(\ref{avlogdist})$.

\section{Acknowledgments}\label{acknow}
We would like to thank Andreas Str\"ombergsson who pointed out Petrov's theorem 
$\ref{Petrov}$ to us.

\providecommand{\href}[2]{#2}
\providecommand{\arxiv}[1]{\href{http://arxiv.org/abs/#1}{arXiv:#1}}
\providecommand{\url}[1]{\texttt{#1}}
\providecommand{\urlprefix}{URL }

\end{document}